\documentclass[10pt]{iopartm}
\usepackage{calrsfs,amsmath,amssymb,graphicx,array,accents,natbib,figlatex,epsfig}
\usepackage{setspace,wrapfig,marvosym, bbm, stmaryrd,centernot, amsthm, color}

\newtheorem{theorem}{Theorem}[section]
\newtheorem{lemma}[theorem]{Lemma}

\newtheorem{rem}{Remark}

\newcommand{\R}{\mathbb{R}}

\newlength{\kaka}

\newcommand{\ahref}[2]{}

\newcommand{\obs}{^{\text{obs}}}

\newcommand{\beq}{\begin{equation}}
\newcommand{\eeq}{\end{equation}}
\newcommand{\lb}{\label}

\newcommand{\bea}{\begin{eqnarray}}
\newcommand{\eea}{\end{eqnarray}}
\newcommand{\bxr}{\begin{array}}
\newcommand{\exr}{\end{array}}

\newcommand\exs{\hspace*{0.4mm}}

\newcommand\nxs{\hspace*{-0.2mm}}

\newcommand{\norms}[1]{\parallel\! #1 \!\parallel}

\newcommand{\bU} {\boldsymbol{U}}
\newcommand{\bSig} {\boldsymbol{\Sigma}}

\newcommand{\bC} {\boldsymbol{C}}
\newcommand{\bK} {\boldsymbol{K}}

\newcommand{\bM} {\boldsymbol{M}}

\newcommand{\bn} {\boldsymbol{n}}
\newcommand{\ba} {\boldsymbol{a}}

\newcommand{\bq} {{\boldsymbol{q}}}
\newcommand{\ff}{^{\text{\tiny f}}}

\newcommand{\bx} {\boldsymbol{x}}

\newcommand{\be} {\boldsymbol{e}}
\newcommand{\bz} {\boldsymbol{z}}
\newcommand{\bg} {{\boldsymbol{g}}}

\newcommand{\bd} {\boldsymbol{d}}
\newcommand{\bI} {\boldsymbol{I}}

\newcommand{\pff}{\boldsymbol{u}^{\textit{i}}}

\newcommand{\sip} {\!\cdot\!}

\newcommand{\bzero}{\boldsymbol{0}}

\newcommand{\bu} {\boldsymbol{u}}
\newcommand{\tbu} {\tilde{\boldsymbol{u}}}
\newcommand{\bt} {{\boldsymbol{t}}}

\newcommand{\bbu} {\llbracket {\tbu}\rrbracket}
\newcommand{\bbut} {\llbracket {\bu}\rrbracket}
\newcommand{\brbu} {\llbracket {\breve\bu}\rrbracket}

\newcommand{\bxi} {\boldsymbol{\xi}}

\def\Xint#1{\mathchoice
{\XXint\displaystyle\textstyle{#1}}%
{\XXint\textstyle\scriptstyle{#1}}%
{\XXint\scriptstyle\scriptscriptstyle{#1}}%
{\XXint\scriptscriptstyle\scriptscriptstyle{#1}}%
\!\int}
\def\XXint#1#2#3{{\setbox0=\hbox{$#1{#2#3}{\int}$}
\vcenter{\hbox{$#2#3$}}\kern-.5\wd0}}

\def\dashint{\Xint-}
\newcommand{\OOd}{\Omega}

\newcommand{\hs}{\hspace{0.01 in}}

\newcommand{\bPhi}{\boldsymbol{\Phi}}

\begin{document}
%=======================================================================================================%

\begin{center}
\title[]{A synoptic approach to the seismic sensing of heterogeneous fractures: from geometric reconstruction to interfacial characterization}
\author{Fatemeh Pourahmadian$^1$, Bojan B. Guzina$^1$ and Houssem Haddar$^2$}
\address{$^1$ Department of Civil, Environmental and Geo-Engineering, University of Minnesota, Twin Cities, USA}
\address{$^2$ INRIA, Ecole Polytechnique (CMAP) and Universit\'e Saclay Ile de France, 91128, Palaiseau, France}
\ead{guzin001@umn.edu}
\end{center} 

\begin{abstract}
A non-iterative waveform sensing approach is proposed toward (i)~geometric reconstruction of penetrable fractures, and (ii)~quantitative identification of their heterogeneous contact condition by seismic i.e.~elastic waves. To this end, the fracture support $\Gamma$ (which may be non-planar and unconnected) is first recovered without prior knowledge of the interfacial condition by way of the recently established approaches to non-iterative waveform tomography of heterogeneous fractures, e.g. the methods of generalized linear sampling and topological sensitivity. Given suitable approximation~$\breve\Gamma$ of the fracture geometry, the jump in the displacement field across~$\breve\Gamma$ i.e.~the fracture opening displacement (FOD) profile is computed from remote sensory data via a regularized inversion of the boundary integral representation mapping the FOD to remote observations of the scattered field. Thus obtained FOD is then used as input for solving the traction boundary integral equation on~$\breve\Gamma$ for the unknown (linearized) contact parameters. In this study, linear and possibly dissipative interactions between the two faces of a fracture are parameterized in terms of a symmetric, complex-valued matrix $\bK(\bxi)$ collecting the normal, shear, and mixed-mode coefficients of specific stiffness. To facilitate the high-fidelity inversion for $\bK(\bxi)$, a 3-step regularization algorithm is devised to minimize the errors stemming from the inexact geometric reconstruction and FOD recovery. The performance of the inverse solution is illustrated by a set of numerical experiments where a cylindrical fracture, endowed with two example patterns of specific stiffness coefficients, is illuminated by plane waves and reconstructed in terms of its geometry and heterogeneous (dissipative) contact condition.  
\end{abstract}


%=======================================================================================================%
\section{Introduction} \label{sec1}
%=======================================================================================================%

Geometric and interfacial properties of fractures in rock and like quasi-brittle solids (e.g.~concrete and composites) are the subject of critical importance to a wide spectrum of scientific and technological facets of our society including energy production from natural gas and geothermal resources~\cite{Verdon2013(2),Verdon2013,Taro2010}, seismology~\cite{McLas2012}, non-destructive evaluation (NDE)~\cite{Doqu2015}, hydrogeology~\cite{Cook1992}, and environmental protection~\cite{Place2014}. One particular quantity embodying the fracture's linearized contact law is the so-called {\emph{specific stiffness matrix}} $\bK$, relating the surface traction to the jump in displacement across the interface. In practical terms, the spatial heterogeneity of the contact parameters reflected in~$\bK(\bxi)$ -- due to e.g.~variable distribution of normal stress~\cite{Martel1989}, may be responsible for progressive failure along discontinuities that may occur well before the frictional resistance of the entire interface is surpassed~\cite{Bobet2005}. This may lead to a catastrophic failure of dams, tunnels and slopes~\cite{Gu1993,Eber2004,Bakun2011}, particularly when fractures show slip-weakening behavior~\cite[e.g.][]{French2014} while the underlying design is based upon \emph{averaged} contact properties.  It is further shown in \cite{Hedayat2014} that the onset of slip along an interface can be identified via temporal variation of the fracture's specific stiffness in shear direction. Accordingly, real-time monitoring of $\bK(\bxi)$'s evolution may not only serve as an early indicator of the interfacial instability and failure, but may also help understand the mechanism of shallow earthquakes~\cite{Jaeg2009}. Active seismic sensing of fractures' contact condition has also come under a spotlight in energy production from unconventional resources~\cite{Verdon2013(2),Verdon2013}, owing to the strong correlation between the hydraulic conductivity of a fracture network and the spatiotmporal variations of $\bK$~\cite{pyrak2016}. For sensing purposes, one should bear in mind that the fracture's response to given activation is driven not only by its contact condition $\bK$, but also by its geometry which is not limited to the planar condition~\cite[e.g.][]{Brac1963,Bart1976,Thom1993}. Thus, the objective of this research is to establish a robust framework for the seismic waveform sensing of heterogeneous fractures, that is capable of resolving both their geometry and interfacial characteristics \emph{without iterations}. Traditionally, seismic waves have been used in the context of acoustic emission (AE)~\cite{Lock1993,Shah1995} to monitor the progression of evolving fractures via the detection of underlying microseismic events, whose energy -- captured by the receivers -- is used to track the failure process. Such ``passive'' sensing approach is, however, ineffective when trying to either image in situ fractures or to assess the fracture's interfacial condition. Another approach, motivating this research, is the concept of \emph{active seismic sensing} applied to fracture identification. This approach, where the discontinuity is ``illuminated'' by an external seismic source, carries the potential of simultaneous fracture imaging and characterization thanks to the sensitivity of scattered waves to the interfacial condition. 

In general the relationship between the wavefield scattered by an obstacle and its geometry and mechanical characteristics is nonlinear, which invites two overt solution strategies: (i) linearization via e.g.~Born approximation and ray theory~\cite{Blei2001}, or ii)~pursuit of the nonlinear minimization approach~\cite{Viri2009}. Over the past two decades, however, a number of \emph{sampling methods} have emerged that consider the nonlinear nature of the inverse problem in an iteration-free way.  In the context of extended scatterers, examples of such paradigm include the linear sampling method (LSM)~\cite{Colt1996,Fiora2003,Mad2007}, the factorization method (FM)~\cite{Kirsch2008, Boukari2013}, the generalized linear sampling method (GLSM)~\cite{Audibert2015,Fatemeh2016, BookCCH2016}, the concept of topological sensitivity (TS)~\cite{Guzi2004,Chik2007,Belli2009} and the subspace migration technique~\cite{Park2015}. Among these, the TS and GLSM approaches have been recently adapted to permit elastic-wave sensing of heterogeneous fractures~\cite{Fatemeh2015,Fatemeh2016}. 

As indicated earlier, there is a mounting interest in medical diagnosis, target recognition, seismology, and energy production~\cite{Fiora2014,Minato2014,Verdon2013,Place2014} to develop hybrid sensing schemes that also reveal the boundary condition of hidden anomalies. In this vein, it is noteworthy that some anomaly-indicator functionals -- such as those featured by the FM~\cite{Kirsch2008}, (G)LSM~\cite{Fiora2010, Fatemeh2016} and TS~\cite{Fatemeh2015} are largely insensitive to the (unknown) boundary condition of a hidden anomaly. For instance, \cite{Fiora2014,Hasse2012} demonstrate that the LSM is successful in reconstructing electromagnetic obstacles and cracks regardless of their boundary condition. Recent advancements on the recovery of boundary (or interfacial) conditions are, on the other hand, mostly \emph{optimization-based} as proposed in the context of acoustic and electromagnetic inverse scattering. Hitherto, a variational method is proposed in~\cite{Colton2004} within the framework of the LSM to determine the essential supremum of the electrical impedance at the boundary of partially coated obstacles. More recently,~\cite{Fiora2014} combined the LSM with iterative algorithms as a tool to expose the surface properties of obstacles from acoustic and electromagnetic data. In elastodynamics, \cite{Minato2014} proposed a Fourier-based approach employing the reverse-time migration and wavefield extrapolation to retrieve the heterogeneous compliance of a planar interface under the premise of high frequency and absence of evanescent waves along the interface.   

In light of the above developments, the twin aim of (non-iterative) seismic imaging and interfacial characterization of heterogeneous fractures is pursued \emph{sequentially} via a 3-step approach where: (1)~the shape reconstruction of penetrable fractures is effected, without prior knowledge of their contact condition, via recently established GLSM~\cite{Fatemeh2016} and TS~\cite{Fatemeh2015} approaches to elastic waveform tomography of discontinuity surfaces; (2) given the reconstructed fracture geometry~$\breve\Gamma$, the fracture opening displacement (FOD) profile is recovered via a double-layer potential representation mapping the FOD to the scattered field observations, and (3)~the specific stiffness profile (as given by its normal, shear, and mixed-mode components) -- is resolved from the knowledge of~FOD and~$\breve\Gamma$ by making use of the traction boundary integral equation written for~$\breve\Gamma$. To help construct a robust inverse solution, a three-step regularization scheme is also devised to minimize the error due to: (i)~compactness of the double-layer potential map used to recover the FOD; (ii)~inexact geometric reconstruction of the fracture surface, and (iii)~presence (if any) of areas on $\breve\Gamma$ with near-zero FOD values.  The performance of the inverse solution is illustrated by a set of numerical experiments assuming seismic illumination in the resonance region. 

%=======================================================================================================%
\section{Preliminaries}\label{MP}
%=======================================================================================================%

\begin{figure}[bp]
\center\includegraphics[width=0.4\linewidth]{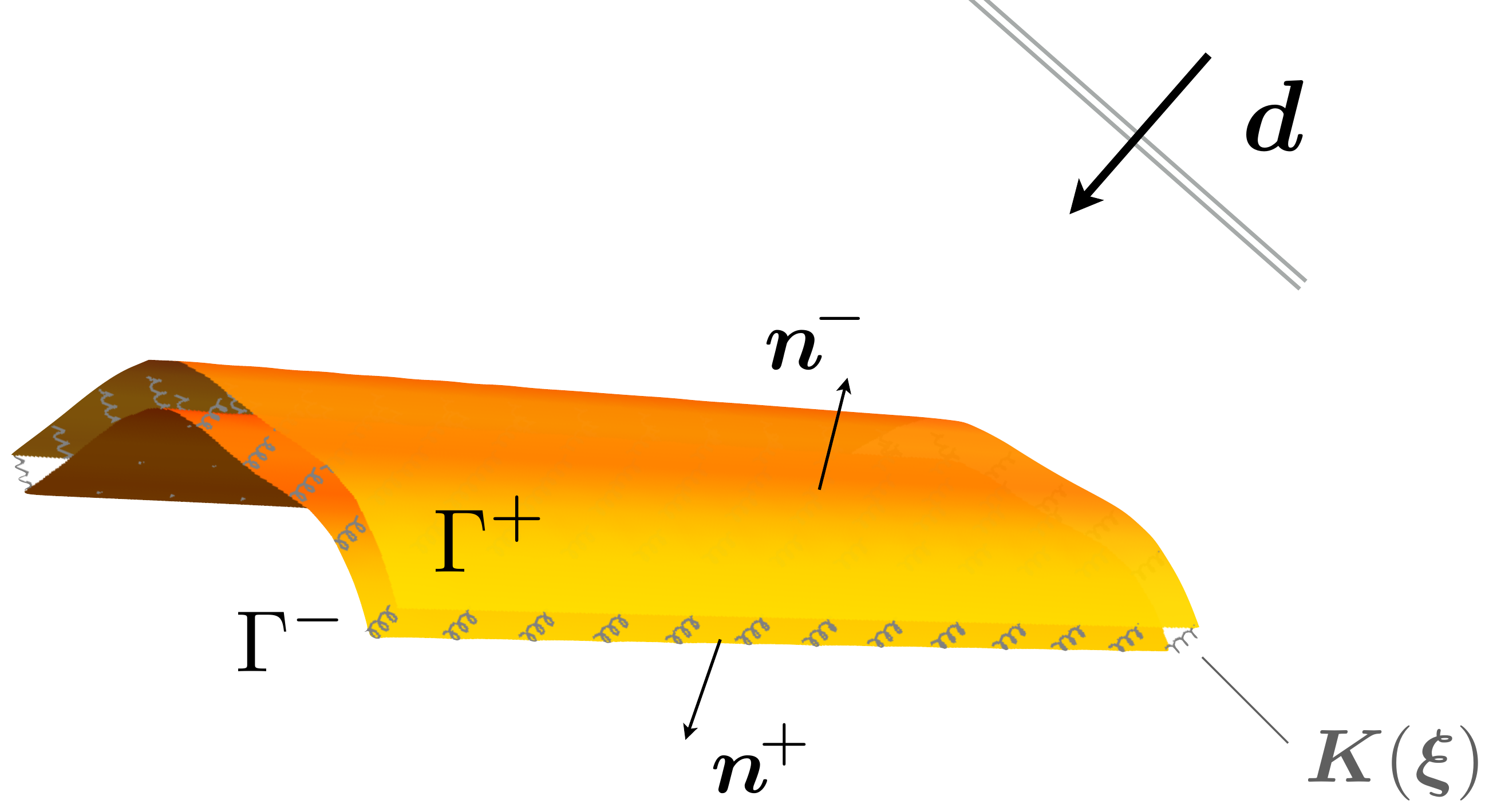} \vspace*{0mm} 
\caption{Fracture~$\Gamma\!\subset\!\mathbb{R}^3$ endowed with a heterogeneous distribution of specific interfacial stiffness $\bK(\bxi)$ is illuminated by a set of plane (compressional and shear) waves propagation in direction $\bd$.} \lb{fig1}
\end{figure} 

With reference to Fig.~\ref{fig1}, consider a heterogeneous (possibly unconnected) fracture $\Gamma\!\subset\!\mathbb{R}^3$ embedded in a homogeneous, isotropic elastic solid endowed with mass density~$\rho$ and Lam\'{e} parameters $\mu$ and~$\lambda$. In the spirit of the linear slip model~\cite{Schoenberg1980} used to describe the seismic response of fractures in rock, the contact condition over $\Gamma$ is given by a $3\times 3$ symmetric matrix, $\bK(\bxi)$, of specific stiffness coefficients -- synthesizing the spatially-varying nature of its rough and possibly multi-phase interface. Assuming time-harmonic seismic illumination, $\bK$ is taken to be complex-valued with~$\Im(\bK) \leqslant \bzero$ in order to allow for energy dissipation at the interface and to ensure the well-posedness of the forward scattering problem~\cite{Fatemeh2016}. 

Let $\Omega$ denote the unit sphere centered at the origin. For a given triplet of vectors $\bd\in\Omega$ and~$\bq_p,\bq_s\!\in\mathbb{R}^3$ such that $\bq_p\!\parallel\bd$ and~$\bq_s\!\perp\!\bd$, consider the case when the fracture is illuminated by a combination of compressional and shear plane waves 
\beq\lb{plwa}
\pff(\bxi) ~=~ \bq_p \exs e^{\textrm{i} k_p \bd\cdot\bxi} \:+\: \bq_s \exs e^{\textrm{i} k_s \bd\cdot\bxi}
\eeq
propagating in direction~$\bd$, where $k_p$ and $k_s=k_p\sqrt{(\lambda\!+\!2\mu)/\mu}$ denote the respective wave numbers. The interaction of~$\Gamma$ with the \emph{incident field}~$\pff$ gives rise to the \emph{scattered field} $\tbu\in H^1_{\mathrm{loc}}(\R^3\backslash\Gamma)^3$ which satisfies 
\beq\lb{GE}
\begin{aligned}
&\nabla \sip (\bC \colon \! \nabla \tbu) \,+\, \rho \exs \omega^2\tbu ~=~ \bzero \quad &\text{in}& \quad {\R^3}\backslash\Gamma, \\*[1mm]
&\bn \cdot \bC \exs \colon \!  \nabla  \tbu~=~ \bK(\bxi) \bbu  \,-\, \bt^{i} \quad &\text{on}& \quad \Gamma. 
\end{aligned}      
\eeq
Here, $\omega^2=k_s^2 \mu/\rho$ is the frequency of excitation; $\bbu=[\tbu^+\!-\tbu^-]$ is the jump in~$\tbu$ across~$\Gamma$, hereon referred to as the fracture opening displacement \textcolor{black}{(FOD)}; 
\beq\label{bC}
\bC \:=\: \lambda\,\bI_2\!\otimes\bI_2 \:+\: 2\mu\,\bI_4  
\eeq
is the fourth-order isotropic elasticity tensor; $\bI_m \,(m\!=\!2,4)$ denotes the $m$th-order symmetric identity tensor; \mbox{$\bt^{i} = \bn \cdot \bC \colon \! \nabla \bu^{i}$} is the incident-field traction vector, and $\bn := \bn^-$ is the unit normal on~$\Gamma$. For clarity it is noted that, owing fo the continuity of the incident field~$\pff$, the jump in the \emph{total field} $\bu=\pff+\tbu$ across~$\Gamma$ equals that of the scattered field, namely  
\beq\lb{equal}
\bbut\;=\;\bbu \qquad\text{on}\quad \Gamma.
\eeq

On uniquely decomposing~$\tbu$ into an irrotational part and a solenoidal part as $\tbu = \tbu_p + \tbu_s$ where 
\beq\label{vpvs}
\tbu_p = \frac{1}{k_s^2\!-\!k_p^2}(\Delta+k_s^2)\tbu, \qquad \tbu_{s} = \frac{1}{k_p^2\!-\!k_s^2}(\Delta+k_p^2)\tbu, 
\eeq
the statement of the forward problem can be completed by enforcing the Kupradze radiation condition 
\beq\lb{KS}
\lim_{r \rightarrow \infty}  \! r \exs 
\Big(\frac{\partial\tbu_p}{\partial r} - \text{i} k_p \tbu_p \Big) = \bzero \quad \mbox{ and } \quad 
\lim_{r \rightarrow \infty}\! r \exs \Big(\frac{\partial \tbu^{s}}{\partial r} - \text{i} k_s \tbu^{s}\Big) = \bzero
\eeq
at infinity, where $r=|\bxi|$ is the distance from the origin. 

%~~~~~~~~~~~~~~~~~~~~~~~~~~~~~~~~~~~~~~~~~~~~~~~~~~~~~~~~~~~~~~~~~~~~~~~~~~~~~~~~~~~~~~~~~~~~~~~~~~~~~~~
\paragraph*{Sensory data.}
%~~~~~~~~~~~~~~~~~~~~~~~~~~~~~~~~~~~~~~~~~~~~~~~~~~~~~~~~~~~~~~~~~~~~~~~~~~~~~~~~~~~~~~~~~~~~~~~~~~~~~~~ 
As shown in~\cite{Martin1993}, any scattered wave $\tbu\in H^1_{\mathrm{loc}}(\R^3 \backslash \Gamma)^3$ solving \eqref{GE}-\eqref{KS} has the asymptotic expansion
\beq\lb{vinf}
\tbu(\bxi) ~=~ \exs \frac{e^{\text{i}k_p r}\!}{4 \pi(\lambda\!+\!2\mu)r} \exs \tbu^\infty_p(\hat\bxi) \:+\: 
\frac{e^{\text{i}k_s r}\!}{4\pi\mu r} \exs \tbu^\infty_s(\hat\bxi) \:+\: O(r^{-2}) \qquad \text{as} \quad r=|\bxi|\to\infty, 
\eeq
where $\hat\bxi=\bxi/r$ is the unit direction of observation, while $\tbu^\infty_p\in L^2(\Omega)^3$ and $\tbu^\infty_s\in L^2(\Omega)^3$ denote respectively the far-field patterns of $\tbu_p$ and $\tbu_s$ that admit~\cite{Fatemeh2016} the integral representation  
\beq\lb{vinf2}
\begin{aligned}
&\tbu_p^\infty(\hat\bxi) ~=~- \text{i} k_p \exs \hat\bxi \int_\Gamma \Big\lbrace \lambda \,  \llbracket \tbu \rrbracket \sip \bn   + 2\mu \big(\bn \sip \hat\bxi \big) \exs \llbracket \tbu \rrbracket \sip \hat\bxi  \exs \Big\rbrace \, e^{-\text{i}k_p \hat\bxi \cdot \bx}  \,\, \text{d}S_{\bx}, \\
& \tbu_s^\infty(\hat\bxi) ~=~ -\text{i}k_s \exs \hat\bxi \exs \times \int_{\Gamma} \Big\lbrace \mu \big( \llbracket \tbu \rrbracket \!\times\! \hat\bxi \exs \big)(\bn \sip \hat\bxi \exs ) \,+\, \mu \big( \bn \!\times \hat\bxi \big) (\llbracket \tbu \rrbracket \sip \hat\bxi) \Big\rbrace \, e^{-\text{i}k_s \hat\bxi \cdot \bx} \,\, \text{d}S_{\bx}.
\end{aligned} 
\eeq
In this setting, we define the far-field pattern of~$\tbu$ by
\beq\lb{far-field}
\tbu^\infty := \tbu^\infty_p \oplus \tbu^\infty_s.
\eeq
For the purposes of seismic fracture sensing, $\tbu^\infty$ is monitored over a subset $\Omega\obs\subseteq\Omega$ of the unit sphere. 

\begin{rem}
The featured framework of elastodynamic scattering in~$\mathbb{R}^3$ is introduced for the reasons of convenience only. In general, the ensuing approach to quantitative reconstruction of the specific stiffness profile~$\bK(\bxi)$, once the fracture geometry has been resolved, applies equally to situations when the reference elastic domain~$D\subset\mathbb{R}^3$ is semi-infinite, bounded, or heterogeneous (e.g. layered half-space). In this vein, the analysis also caters for sensing configurations entailing near-field observations of the scattered field~$\tbu$ over a limited-aperture surface~$S\obs\subset\overline{D}$. 
\end{rem}

%=======================================================================================================%
\section{Elastic-wave sensing of heterogeneous fractures}
%=======================================================================================================%

\begin{figure}[!h]
\vspace*{1mm}
\begin{center}
\includegraphics[width=0.98\linewidth]{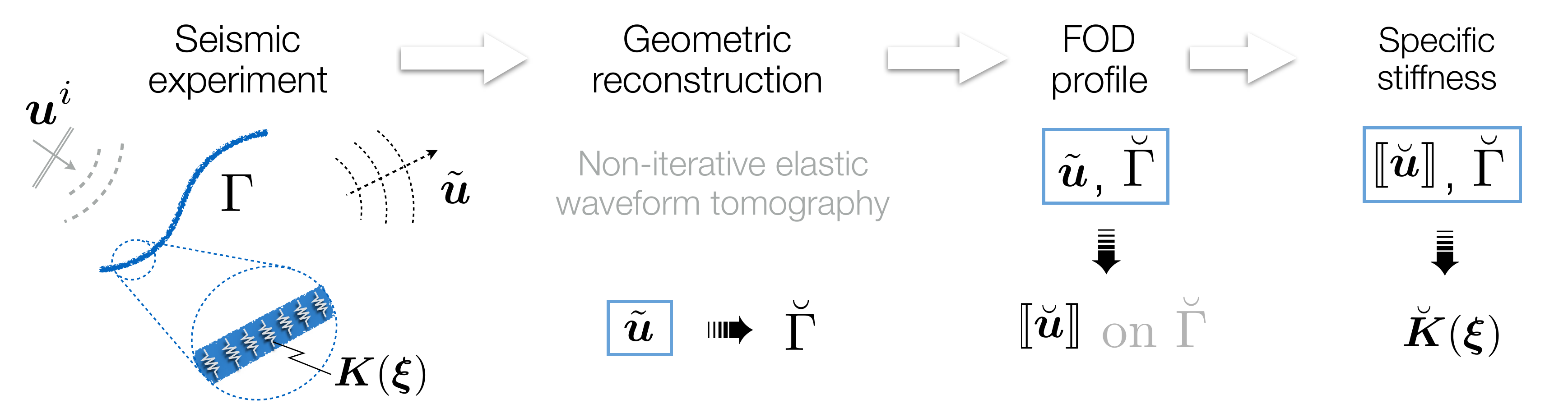}
\end{center} \vspace*{-4mm}
\caption{\small{Three-step approach to non-iterative geometric reconstruction and interfacial characterization of heterogeneous fractures by elastic waves.}}
\label{inversion3}\vspace*{-4mm}
\end{figure}

In what follows, the model problem in Section~\ref{MP} is used as a basis to describe the inverse solution schematically shown in Fig.~\ref{inversion3}, where the sensory waveform data (in this case the far-field pattern $\tbu^\infty$) provide an input to the 3-step approach for geometric reconstruction and interfacial characterization of subsurface discontinuities in an elastic solid, e.g.~natural or hydraulic fractures in rock. In this framework,
\begin{itemize}
\item The \emph{fracture geometry}~$\breve\Gamma$ is reconstructed, irrespective of its contact condition, via a suitable approach to non-iterative seismic waveform tomography~\cite{Fatemeh2015, Fatemeh2016}; 

\item The \emph{FOD profile}, $\brbu$, over the reconstructed fracture support~$\breve\Gamma$ is recovered from the germane boundary integral representation (double-layer potential) of the scattered field, and 

\item The \emph{specific stiffness profile} $\breve\bK(\bxi)$ -- as given by its normal, shear, and mixed components -- is recovered from the knowledge of~$\brbu$ and~$\breve\Gamma$.  
\end{itemize}
These basic steps are elucidated in the sequel. 

%--------------------------------------------------------------------------------------------------%
\subsection{Geometric fracture reconstruction} \lb{TSP} 
%--------------------------------------------------------------------------------------------------%

\noindent The essential first ingredient toward comprehensive sensing of heterogeneous fractures is an ability to decipher the observed elastic waveforms toward geometric reconstruction of germane discontinuity surfaces. By building on the earlier works in scalar inverse scattering~\cite{Boukari2013,Audibert2014}, electromagnetism~\cite{park2012}, and elastic-wave sensing of impenetrable (traction-free) discontinuities~\cite{Belli2009,belli2013}, recent research efforts have shown a path toward non-iterative elastic waveform tomography of \emph{penetrable fractures} (irrespective of their contact condition) via the TS approach~\cite{Fatemeh2015} -- as corroborated by high-frequency simulations and laboratory observations~\cite{Tokm2013,Fatemeh2015(2)}, and the GLSM paradigm~\cite{Fatemeh2016}. To provide an example and to help establish an explicit setting for the discussion, the GLSM approach to elastic-wave imaging of heterogeneous fractures is summarized next. 
% \begin{rem}
% Both GLSM and TS techniques are generic in that they apply to finite, semi-infinite, or infinite (homogeneous or piecewise-homogeneous) background domains~$\Omega$ -- whose elastic and density properties are known beforehand. 
% \end{rem}

\subsubsection{Generalized Linear Sampling Method.}\lb{GLSM}  
For a given vector density $\bg=\bg_p\oplus\bg_s \in L^2(\OOd)^3$ comprised of its compressional ($\bg_p$) and shear ($\bg_s$) wave components mimicking the decomposition of~$\bq$ in~(\ref{plwa}), consider the elastic Herglotz wave function~\cite{Dassios1995} 
\beq\lb{HW}
\pff_\bg(\bxi) ~: =~  \int_{\OOd} \bg_p(\bd) \exs  e^{\textrm{i} k_p \bd \cdot \bxi} \,\, \text{d}S_{\bd} ~\; \oplus \, \int_{\OOd} \bg_s(\bd) \exs  e^{\textrm{i} k_s \bd \cdot \bxi} \,\, \text{d}S_{\bd},  \qquad \bxi \in \R^3, 
\eeq
and define the \emph{far-field operator} $F: L^2(\OOd)^3 \to L^2(\OOd)^3$ by
\beq\lb{ffo0}
F(\bg) \;:=~ \textcolor{black}{\tilde\bu_{\bg}^\infty,} 
\eeq
where~$\tilde\bu_{\bg}^\infty$ is the far-field pattern~\eqref{far-field} of~$\tilde\bu\!\in\!H^1_{\mathrm{loc}}(\R^3 \backslash \Gamma)^3$ solving~\eqref{GE} and~\eqref{KS} with data $\bt^{i} = \bn \cdot \bC \colon \! \nabla \bu_\bg^{i}$ on~$\Gamma$.  
In this setting, the idea of the GLSM is to construct a nearby solution of the (ill-posed) far-field equation 
\beq\lb{LSM}
F \bg ~\simeq~ \bPhi_L^\infty,
\eeq
where $\bPhi_L^\infty$ is the far-field pattern of a \emph{test} radiating field, generated by a small \emph{trial} fracture~$L\subset\mathbb{R}^3$ with prescribed FOD~$\ba\in \tilde{H}^{1/2}(L)^3$ (see~\cite{Fatemeh2016} for details). Without loss of generality, $L$ can be taken as a vanishing penny-shaped fracture at~$\bz\in\mathbb{R}^3$ with normal~$\bn\in\Omega$ and constant (mode~I) FOD profile $\ba\propto\bn$, in which case~\eqref{vinf2} yields  
\beq\lb{Phi-inf-num}
\bPhi_L^\infty(\hat\bxi) ~=~ -  \Big(\text{i} k_p \,  \hat\bxi \exs  \big[ \exs  \lambda+2\mu \exs (\textrm{\bf{n}} \cdot \hat\bxi)^2  \exs \big] \exs e^{-\text{i}k_p \hat\bxi \cdot \bz} \;\oplus\;
 2 \text{i} \mu \exs k_s \,\hat\bxi \times \nxs (\textrm{\bf{n}} \times\hat\bxi)\exs   (\textrm{\bf{n}}\cdot\hat\bxi) \, e^{-\textrm{\emph{i}} k_s \hat\bxi \cdot \bz}  \Big).
\eeq
With reference to~\eqref{LSM}, a heterogeneous fracture~$\Gamma$ with arbitrary (linear) contact condition can then be reconstructed~\cite{Fatemeh2016} from the support of the \emph{GLSM characteristic function}
\beq 
\mathcal{C}\big(L(\bz,\bn)\big) ~=~
\Big(\big\|(F_\sharp)^{\frac{1}{2}} \exs \bg_{\alpha,\delta}^L\big\|^2  \exs+\, 
\delta \big\|\bg_{\alpha,\delta}^L\big\|^2\Big)^{-1/2}, \qquad \bz\in\mathbb{R}^3, ~~ \bn\in\Omega
\eeq   
where 
\beq\lb{Fs}
F_\sharp \,:=\; \tfrac{1}{2}|(F+F^*)| \:+\: \tfrac{1}{2\textrm{i}}(F-F^*) 
\eeq  
is a self-adjoint surrogate for~$F$~\cite{Kirsch2008}, and~$\bg_{\alpha,\delta}^L$ minimizes the penalized least-squares functional
\beq  \lb{RJ-alph}
J_\alpha^\delta(\exs\bg;\exs\bPhi_L^\infty) ~ \colon \!\!\! =~ \!   \norms{\nxs F\bg\,-\,\bPhi_L^\infty \nxs}^2 + \,\, \alpha \exs \big(\!\norms{\nxs(F_\sharp)^{\frac{1}{2}} \exs \bg \nxs}^2 +\,\,\exs \delta  \!  \norms{\nxs \bg \nxs}^2 \! \big), \qquad \bg\in L^2(\OOd)^3
\eeq
that features an absolute measure, $\delta>0$, of perturbation in the data (given by the far-field operator~$F$) and a penalization parameter~$\alpha=\alpha(\delta)>0$ as in~\cite{Audibert2014}. 
  
\begin{rem}
For generality, it should be noted that geometric fracture reconstruction can also be accomplished by other application-specific inversion schemes, e.g.~microseismic imaging --  which deploys travel time inversion or migration of seismic events generated during the fracturing process~\cite{Maxw2014}. The featured TS and GLSM techniques are, however, part of a  more general waveform inversion platform which is (i)~computationally efficient due to its non-iterative nature; (ii)~self-contained for it uses the common set of sensory data for both geometric reconstruction and interfacial characterization, and (iii)~robust by providing significant flexibility in terms of the sensing configuration and requiring no a~priori knowledge on the fracture geometry nor its contact condition. 
\end{rem}

%--------------------------------------------------------------------------------------------------%
\subsection{Inversion of the FOD profile}\lb{FOD}
%--------------------------------------------------------------------------------------------------%

\noindent Given a geometric reconstruction of the fracture surface $\breve\Gamma$ -- obtained as described in Section~\ref{TSP}, a double-layer (elastodynamic) potential representation of the scattered field~\cite{Bon1999} serves as \emph{a map} $\breve{\bM}\!\!:\tilde{H}^{1/2}(\breve\Gamma)^3\to L^2(\Omega\obs)^3$ relating the sought FOD, $\brbu$, to the sensory data $\tbu^\infty$, namely 
\beq
\lb{FOD1}
\breve{\bM} \brbu(\hat\bxi) ~=~ \tbu^\infty (\hat\bxi), \qquad  
\breve{\bM}\brbu(\hat\bxi) ~:=\, \int_{\breve\Gamma} \, \bSig^\infty(\hat\bxi,\bx) \colon \! \big\{\brbu(\bx) \otimes \breve\bn(\bx)\big\} \,\, \text{d}S_{\bx},  
\qquad \hat\bxi\in\Omega\obs \vspace*{-1mm}
\eeq    
where $\breve\bn = \breve\bn^-$ is the unit normal on~$\breve\Gamma$, and $\bSig^\infty$ is the \emph{far-field pattern} as $|\bxi|\to\infty$ of the elastodynamic fundamental stress tensor $\bSig(\bxi,\bx)$~\cite{Ach2003}. In terms of dyadic notation and Einstein summation convention, the latter quantity can be written as 
\beq\label{sfund2}
\bSig(\bxi,\bx) ~=~ \Sigma_{ij}^\ell(\bxi,\bx)\, \be_\ell\otimes\be_i\otimes\be_j, \qquad \bxi,\bx\in\mathbb{R}^3, \quad \bxi\ne\bx, \quad i,j,\ell\in\{1,2,3\}
\eeq
where~$\Sigma_{ij}^\ell\!=\Sigma_{ji}^\ell$ signify the components of the Cauchy stress tensor at $\bxi$ due to point force $\be_\ell$ (i.e.~the unit vector in the $\ell$th coordinate direction) acting at $\bx$. For completeness, it is noted (roughly speaking) that the Sobolev space $\tilde{H}^{1/2}(\breve\Gamma)$ denotes the space of functions in ${H}^{1/2}(\breve\Gamma)$ whose extension by zero to a larger Lipshitz surface, $\Lambda\supset\breve\Gamma$, gives functions that belong to ${H}^{1/2}(\Lambda)$ \cite{Fatemeh2016}.  The space $\tilde{H}^{1/2}(\breve\Gamma)$ is also the dual of  $H^{-1/2}(\breve\Gamma)$.

\begin{lemma}\lb{compact} 
Linear operator $\breve{\bM}\!\!:\tilde{H}^{1/2}(\breve\Gamma)^3\to L^2(\Omega\obs)^3$ introduced in~(\ref{FOD1}) is compact and injective (see also Lemma~5.2 in~\cite{Fatemeh2016}).
\end{lemma}
\begin{proof}
Thanks to the explicit asymptotic expansion of $\bSig$ as $|\bxi|\to\infty$~\cite[e.g.][]{Ach2003} and~(\ref{vinf}), one finds that 
\begin{multline}\lb{siginf}
\Sigma_{ij}^{\ell,\infty}(\hat\bxi,\bx) \;=\; \text{i}k_p (2\mu\exs\hat\xi_i\exs \hat\xi_j + \lambda \delta_{ij})  \hat\xi_\ell \, e^{-\text{i}k_p\hat\bxi\cdot\bx} 
\,\oplus \, \text{i}k_s \mu (\delta_{i\ell}\exs\hat\xi_j + \delta_{j\ell}\exs\hat\xi_i - 2 \hat\xi_i \exs\hat\xi_j\exs \hat\xi_\ell)  \, e^{-\text{i}k_s\hat\bxi\cdot\bx}, \\ \quad i,j,\ell\in\{1,2,3\}
\end{multline}
which are holomorphic over~$\Omega\times\mathbb{R}^3$. As a result, integral operator~$\breve{\bM}$ has a smooth kernel and is therefore compact from~$\tilde{H}^{1/2}(\breve\Gamma)^3$ to~$L^2(\Omega\obs)^3$. To establish the injectivity of~$\breve{\bM}$, one may observe as in~\cite{Fatemeh2016} that the vanishing far-field pattern~\eqref{vinf2} requires that~$\tbu\!=\!\bzero$ in~$\mathbb{R}^3\backslash\breve\Gamma$ by the Rellich Lemma, and consequently that~$\brbu\!=\!\bzero$ thanks to the unique continuation principle and the fundamental property of double-layer potentials by which~$\brbu=\llbracket\tbu\rrbracket$ on~$\breve\Gamma$. 
\end{proof}
\begin{rem} In the sequel, it is assumed that $\breve{\bM}\!\!:\tilde{H}^{1/2}(\breve\Gamma)^3\to L^2(\Omega\obs)^3$ has a dense range, a hypothesis that is supported by Lemma~5.3 in~\cite{Fatemeh2016} giving sufficient conditions on the fracture geometry~$\breve\Gamma$ and excitation frequency~$\omega$ for the range denseness of~\eqref{FOD1}.
\end{rem}
\begin{rem}\label{int-waves}In physical terms, the compactness of~$\breve{\bM}$ is reflected by the presence of interface waves~\cite{Pyr1987} on~$\breve\Gamma$. In theory these local waves, propagating along the surfaces of discontinuity in an elastic medium, are characterized by an exponential decay~\cite{Stone1924,Pyr1987} with normal distance to the interface. Despite the apparent leakage of such surface-wave energy into the exterior due to the curvature of~$\breve\Gamma$ (if any) and the interaction of interfacial waves with the fracture edge~$\partial\breve\Gamma$, these local FOD modes may have only a marginal fingerprint in the remote sensory data that warrants a custom treatment in the inversion scheme.
\end{rem}

\begin{rem}
For a generic sensing configuration entailing (i) semi-infinite or bounded reference domain~$D\subset\mathbb{R}^3$ and (ii) near-field measurements of the scattered field, $\tbu^\infty$, $\Omega\obs$, and $\bSig^\infty$ in~(\ref{FOD1}) are suitably replaced by $\tbu$, $S\obs\!\subset\overline{D}$, and $\bSig_D$ -- the germane elastodynamic Green's function. Assuming a physical separation between the fracture surface and the sensing grid, the relevant double-layer potential map $\bM_{\!D}$ can likewise be shown (on a case-by-case basis) to be compact thanks to the smoothness of~$\bSig_D$.  
\end{rem}

\noindent \emph{Discretization.} Taking~$\Omega\obs$ as a union of discrete sensing directions -- as may be the case in a physical experiment, and allowing the reconstructed fracture surface $\breve\Gamma$ to be arbitrarily complex, a discrete version of $\breve{\bM}\brbu$ may be obtained via suitable discretization of $\breve\Gamma$ and FOD in terms of surface i.e.~boundary elements~\cite{Bon1999,Fatemeh2015}. As a result,~(\ref{FOD1}) can be recast via the collocation method as     
\beq
\lb{FOD2}
\breve{\text{\bf M}} \exs \llbracket \breve{\text{\bf u}} \rrbracket ~=~ \tilde{\text{\bf u}}^\infty, 
\eeq  
where~$\breve{\text{\bf M}}$ is a $3N^\text{obs}\times 3N^\text{nds}$ coefficient matrix; $\llbracket \breve{\text{\bf u}}\rrbracket = \llbracket \breve{\tilde{\text{\bf u}}}\rrbracket$ is a vector sampling the FOD at $\,N^\text{nds}$ geometric nodes over~$\breve\Gamma$ (see~\cite[Fig.~11]{Fatemeh2016}), and $\tilde{\text{\bf u}}^\infty$ collects the far-field pattern~(\ref{far-field}) of the scattered field in~$N^{\text{obs}}$ sensing directions over~$\Omega\obs$. % Note that the components of~$\breve{\text{\bf M}}$ are evaluated without regularization thanks to the analyticity of~$\bSig^\infty$ in~(\ref{FOD1}).
\vspace*{2mm}

\noindent \emph{Solution.}~In this setting, the idea is to have $N^\text{nds}<N^{\text{obs}}$ and to compute the FOD profile on $\breve{\Gamma}$ from the sensory data~$\tilde{\text{\bf u}}^\infty$ by solving an overdetermined linear system~(\ref{FOD2}). As a result, the \emph{spatial resolution} of the FOD reconstruction will be inherently limited (in addition to other factors) by the number of observation directions on~$S\obs$. Also note that~(\ref{FOD2}) can be solved anew for each available direction $\bd_p\in\Omega$ ($p=1,\dots P$) of the incident plane wavefield~(\ref{plwa}), where $\breve{\text{\bf M}}$ is invariant i.e. source-independent. \vspace*{2mm}

\noindent \emph{Regularization.}~A critical point in recovering the FOD -- which subsequently affects the inversion of the specific stiffness -- is that~(\ref{FOD1}) is \emph{ill-posed} due to the compactness of~$\breve{\bM}$. In the context of discrete statement~(\ref{FOD2}), $\breve{\text{\bf M}}$ may accordingly contain unacceptably small singular values, whose number will (for a given fracture geometry and contact characteristics) depend on the properties of the incident field~$\bu^i$. To provide a physical interpretation of such behavior, consider the singular value decomposition (SVD) of~$\breve{\text{\bf M}}$ as
\beq\lb{SVDM}
\breve{\text{\bf M}} \;=\; \text{\bf U}_{\text{\tiny M}} \boldsymbol{\Lambda}_{\text{\tiny M}} \text{\bf V}_{\!\text{\tiny M}}^*, \qquad 
\eeq
where $\text{\bf U}_{\text{\tiny M}}$ (resp. $\text{\bf V}_{\!\text{\tiny M}}$) collects the left (resp. right) eigenvectors of~$\breve{\text{\bf M}}$, and $\boldsymbol{\Lambda}_{\text{\tiny M}}$ contains its singular values. In this setting, the \emph{right eigenvectors} $\text{\bf V}_{\!\text{\tiny M},q}$  ($q\!=\!\overline{1,Q}$) corresponding to the~$Q$ smallest (unacceptably small) singular values of~$\breve{\text{\bf M}}$, $\Lambda_q\leqslant\epsilon$ (which reflect the compactness of~$\breve\bM$), can be affiliated with the \emph{interface wave modes}~\cite{Pyr1987} on~$\breve\Gamma$, see Remark~\ref{int-waves}. In the context of~(\ref{FOD2}), the ill-posedness of~(\ref{FOD1}) can be dealt with in a standard way via e.g. truncated SVD or Tikhonov regularization aided by the Morozov discrepancy principle~\cite{Kirc2011}(which takes the penalty parameter to be commensurate with the level of noise in the data). It will be shown later, however, that the approximating effect of such regularization can be mitigated assuming the availability of sensory data ($\tilde{\text{\bf u}}^\infty$) for \emph{multiple} incident fields -- which is a customary premise for most inverse scattering solutions. 

% This problem may arise due to either a \emph{limited ``viewing'' aperture} furnished by~$S\obs$, or the emergence of~\emph{interfacial scattered waves} -- propagating along the fracture surface~\cite{Pyr1987} -- that cannot be sensed on~$S\obs$ (see Section~\ref{Numexpp}). Accordingly, (\ref{FOD2}) can be solved via suitable regularization, e.g.~Tikhonov regularization aided by the Morozov discrepancy principle (which takes the regularization parameter to be commensurate with the level of noise in the data)~\cite{Kirc2011}. As examined next, the span of eigenvectors corresponding to unacceptably small singular eigenvalues of~$\breve{\text{\bf M}}$ (eliminated by regularization) determines the FOD subspace \emph{to be avoided} when reconstructing the profile of specific stiffnesses.

%--------------------------------------------------------------------------------------------------%
\subsection{Inversion of the specific stiffness profile} \lb{sstif}
%--------------------------------------------------------------------------------------------------%

\noindent The last step in the proposed inverse scheme is to substitute the identified FOD, $\brbu$, into the fracture's boundary condition on~$\breve\Gamma$, and to solve thus-obtained equation for the specific stiffness profile $\breve\bK(\bxi)$. In this vein, the true contact condition~(\ref{GE}) on $\Gamma$ is recast over the reconstructed fracture surface~$\breve\Gamma$ as
\beq \lb{BCFRC}
\breve\bK(\bxi) \exs \hs \brbu ~=~ \breve{\tilde{\bt}} \,+\, \exs \breve{\bt}^i,   \qquad \bxi \in \breve\Gamma, 
\eeq
where $\breve{\bt}^i \!=\! \breve\bn\sip\bC \! \exs \colon \!\! \nabla\bu^i \exs$ is the incident-field traction on $\breve\Gamma$, while $\breve{\tilde{\bt}}$ is the scattered-field traction on $\breve\Gamma$ expressed in terms of $\llbracket \breve{\tbu}\rrbracket=\brbu$ by invoking the traction boundary integral equation (TBIE)~\cite{Fatemeh2015, Bon1999} as a map $\boldsymbol{T}\!\!:\tilde{H}^{1/2}(\breve\Gamma)\to H^{-1/2}(\breve\Gamma)$ such that 
\beq
\lb{TBIE}  
\begin{aligned}
\breve{\tilde{\bt}}(\bxi) ~=~ \breve{\boldsymbol{T}}\brbu \;:=~ - & \exs \breve\bn(\bxi)  \cdot  \bC  :  \dashint_{\breve\Gamma} \, \bSig(\bxi,\bx)  : \nxs \boldsymbol{D}_{\!\bx} \brbu (\bx) \, \text{d}S_{\bx} 
~+~  \\[1 mm]
& \qquad \rho \exs \omega^2 \exs \breve\bn(\bxi)  \cdot  \bC \colon \!\!\! \int_{\breve\Gamma} \, \bU(\bxi,\bx) \cdot \big( \brbu \nxs \otimes \nxs \breve\bn \big)(\bx) \, \text{d}S_{\bx},  \qquad    \bxi \in \breve\Gamma, 
\end{aligned}
\eeq 
where $\bU=U_{i}^\ell(\bxi,\bx)\, \be_i\otimes\be_\ell$ and~$\bSig=\Sigma_{ij}^\ell(\bxi,\bx)\, \be_i\otimes\be_j\otimes\be_\ell$ denote respectively the elastodynamic displacement and stress fundamental solution at~$\bxi\in\mathbb{R}^3$ due to point force acting at $\bx\in\mathbb{R}^3$ (see~\cite[Appendix A]{Fatemeh2015},~\cite{Guzin2001,Bon1999}); $\,\dashint$~signifies the Cauchy-principal-value integral, and $\boldsymbol{D}_{\!\bx}$ is the tangential differential operator given by
\beq\lb{dod}
\boldsymbol{D}_{\!\bx}(\boldsymbol{f}) = D_{kl} (f_m)   \, \be_l\otimes\be_m\otimes\be_k, \qquad 
D_{kl}(f_m)  = \breve{n}_k f_{m,l} - \breve{n}_l f_{m,k}, 
\eeq
with $\breve{n}_k=\breve{n}_k(\bx)$ and~$f_{m,k}=\partial f_m/\partial x_k$ in the global coordinate frame $\lbrace \be_1,\be_2,\be_3 \rbrace$.   

\subsubsection{Indirect solution approach.} Given $\brbu$ on $\breve\Gamma$, the scattered-field traction $\breve{\tilde{\bt}}(\bxi)$ in~(\ref{BCFRC}) can be computed without any reference to~$\breve\bK(\bxi)$. This is accomplished by imposing $\llbracket \breve{\tbu}\rrbracket=\brbu$ as a Dirichlet boundary condition on~$\breve\Gamma$, and solving the resulting (exterior) boundary value problem in the reference domain ($\mathbb{R}^3$ in the present study, or $D\supset\breve\Gamma$ in a more general case). With the knowledge of $\breve{\tilde{\bt}}$ and $\breve{\bt}^i$ on~$\breve\Gamma$ at hand, $\breve{\bK}(\bxi)$ can then in principle be solved from~(\ref{BCFRC}). In particular: 

\begin{itemize}
\item If $\breve\bK$ is assumed to be \emph{diagonal} in the fracture's local coordinates, namely $\breve\bK\!:=\!\text{diag}(\breve{\kappa}_n,$ $\breve{\kappa}_{s_1},\breve{\kappa}_{s_2})$, one arrives via~(\ref{BCFRC}) at the uncoupled system
\[
\breve{\kappa}_p(\bxi) \, \llbracket \breve{u}_p\rrbracket(\bxi) \;=\; \breve{\tilde{t}}_p(\bxi) \;+\; \breve{t}^i_p(\bxi), \qquad p\in\{n,s_1,s_2\}, \quad \bxi\in\breve\Gamma
\]
for the diagonal entries of $\breve{\bK}(\bxi)$ (no summation implied), which is solvable from the sensory data for \emph{a single} incident field -- provided that $\llbracket \breve{u}_\ell(\bxi)\rrbracket$ does not vanish at~$\bxi$. 

\item In situations where $\breve\bK$ is taken to be \emph{fully populated} - implying dilatant contact behavior (i.e. coupling between the normal and shear resistance), one obtains the coupled system
\[
\breve{K}_{pq}(\bxi) \, \llbracket \breve{u}_q\rrbracket(\bxi) \;=\; \breve{\tilde{t}}_p(\bxi) \;+\; \breve{t}^i_p(\bxi), \qquad p,q\in\{n,s_1,s_2\}, \quad \bxi\in\breve\Gamma
\]
which makes use of the Einstein summation convention. In this case there are (recalling the symmetry of~$\breve{\bK}$) six unknowns at given $\bxi\!\in\!\breve\Gamma$, requiring the availability of sensory data for \emph{at least two} incident fields. 
\end{itemize}
This approach has an advantage that it does not require the Green's function for the reference domain, which is convenient when dealing with fracture sensing inside finite bodies, or using numerical solution techniques such as finite element or finite difference methods. 

\subsubsection{Direct solution approach.} One apparent drawback of the foregoing scheme is that it requires an additional forward solution for each incident field, as required to compute the scattered-field traction~$\breve{\tilde{\bt}}(\bxi)$. In situations where the Green's function for the reference domain is available (as in the present case), this problem can be partly circumvented by combining~(\ref{BCFRC}) and~(\ref{TBIE}) as
\beq\lb{dsa1}
\breve\bK(\bxi) \exs \hs \brbu ~=~ \breve{\boldsymbol{T}}\brbu \,+\, \exs \breve{\bt}^i,   \qquad \bxi \in \breve\Gamma,
\eeq
where the right-hand side can be computed as follows. \vspace*{2mm} 

\noindent \emph{Discretization.}~By deploying the collocation method as examined in Section~\ref{FOD}, a discretized version of the integro-differential operator~$\breve{\boldsymbol{T}}\brbu$ in~(\ref{TBIE}) can be computed as $\breve{\text{\bf T}} \llbracket \breve{\text{\bf u}} \rrbracket$, where~$\breve{\text{\bf T}}$ is a $3N^\text{col}\times 3N^\text{nds}$ coefficient matrix; $N^\text{col}$ is the number of points on~$\breve{\Gamma}$ where~(\ref{dsa1}) is collocated, and~$\llbracket\breve{\text{\bf u}}\rrbracket$ is a vector sampling the FOD at $\,N^\text{nds}$ geometric nodes over~$\breve\Gamma$ as before. As examined in the sequel, $N^\text{col}$ -- carrying the parametrization of~$\breve{\bK}(\bxi)$ -- can be either smaller or larger than $N^\text{nds}$. To meet the computational $C^1$ smoothness requirement on FOD in~(\ref{TBIE}) (due to appearance of the tangential differential operator $\boldsymbol{D}_{\!\bx}$) while allowing $\brbu$ on~$\breve\Gamma$ to be parametrized via standard $C^0$ (boundary element) interpolation, the collocation points for solving~(\ref{dsa1}) are conveniently placed in the \emph{interior} of the boundary elements used to represent the geometry and kinematics of the reconstructed fracture surface (see~\cite{Fatemeh2015}, Appendix~B for details). However, since the collocation points belong to the fracture surface, the first integral on the right-hand-side of~(\ref{TBIE}) remains singular and must be properly regularized~\cite{Bon1999,Fatemeh2015}. In this setting, a discretized statement of the contact condition~(\ref{dsa1}) reads 
\beq \lb{Claw}
\breve{\text{\bf K}} \exs \llbracket \breve{\text{\bf u}} \rrbracket ~=~   \breve{\text{\bf T}} \, \llbracket \breve{\text{\bf u}} \rrbracket \,+\, \breve{\text{\bf t}}^i,  
\eeq
where~$\breve{\text{\bf K}}$ is a $3N^\text{col}\times 3N^\text{nds}$ block-diagonal matrix of specific stiffness coefficients, and $\breve{\text{\bf t}}^i$ is the free-field traction vector sampled at $\,N^\text{col}$ collocation points over~$\breve\Gamma$. \vspace*{2mm}

\noindent \emph{Solution.}~To solve~(\ref{Claw}) for the entries of $\breve{\text{\bf K}}$, one may consider the following situations.
\begin{itemize}
\item Assuming $\breve\bK\!:=\!\text{diag}(\breve{\kappa}_n,$ $\breve{\kappa}_{s_1},\breve{\kappa}_{s_2})$, $\breve{\text{\bf K}}$ in~(\ref{Claw}) becomes diagonal and the specific stiffness profile can be evaluated directly at collocation points. In particular, the left-hand side of~(\ref{Claw}) can be recast as  
\beq \lb{Claw2}
\breve{\text{\bf K}} \exs \llbracket \breve{\text{\bf u}} \rrbracket \;:=~   \breve{\text{\bf A}}_{\llbracket \breve{\text{\bf u}} \rrbracket} \exs \breve{\text{\bf k}}, 
\eeq
where~$\breve{\bf A}_{\llbracket\breve{\text{\bf u}}\rrbracket}$ is a $3N^\text{col}\times 3N^\text{col}$ diagonal coefficient matrix given by the reconstructed FOD profile, and $\breve{\text{\bf k}}$ is a $3N^\text{col}$-vector sampling the (normal and shear) specific stiffness coefficients, namely~$\lbrace \breve{\kappa}_n,\breve{\kappa}_{s_1},\breve{\kappa}_{s_2}\rbrace$ at $\,N^\text{col}$ collocation nodes over~$\breve\Gamma$. Assuming $m$ internal collocation points per boundary element in the evaluation of~(\ref{Claw}), one accordingly has $N^\text{col} = m N^\text{el} \gtreqless N^\text{nds}$, where $N^\text{el}$ denotes the number of boundary elements. In situations when~$\breve{\bf A}_{\llbracket\breve{\text{\bf u}}\rrbracket}$ is free of zero or near-zero diagonal elements (the contrary case will be discussed shortly), the substitution of~(\ref{Claw2}) into~(\ref{Claw}) immediately yields a stable solution for the specific stiffness profile, $\breve{\text{\bf k}}$, on~$\breve\Gamma$.   

\item In situations where $\breve\bK$ is taken as fully populated whereby the number of unknowns per contact point is six, the specific stiffness distribution is parametrized using the previously established \emph{geometric interpolation} -- so that the unknowns are evaluated at $N^\text{nds}$ geometric nodes. In this setting, one has
\beq \lb{Claw3}
\breve{\text{\bf K}} \exs \llbracket \breve{\text{\bf u}} \rrbracket \;:=~   \breve{\text{\bf B}}_{\llbracket \breve{\text{\bf u}} \rrbracket} \exs \breve{\text{\bf k}}, 
\eeq   
where $\breve{\text{\bf B}}_{\llbracket \breve{\text{\bf u}} \rrbracket}$ is a $3N^\text{col}\times 6N^\text{nds}$ coefficient matrix given by the recovered FOD profile, and $\breve{\text{\bf k}}$ is a $6N^\text{nds}$-vector collecting the six entries of the symmetric stiffness matrix~$\breve{\text{\bf K}}$ sampled at $\,N^\text{nds}$ geometric nodes over~$\breve\Gamma$. Accordingly, by taking $m$ internal collocation points per element so that $N^\text{col}= m N^\text{el} \geqslant 2N^\text{nds}$, one obtains the sufficient number of equations to solve~(\ref{Claw}) for $\breve{\text{\bf k}}$. Similar to the previous case, this yields a stable solution provided that~$\breve{\text{\bf B}}_{\llbracket \breve{\text{\bf u}} \rrbracket}$ is free of zero or near-zero singular values, which physically implies that no geometric nodes be affiliated with permanent contact (zero FOD) on~$\breve\Gamma$. For completeness, this and other likely sources of ill-posedness jeopardizing the evaluation~$\breve\bK(\bxi)$ are addressed next. 
\end{itemize} 

\subsection{Regularization}

To help construct a robust solution to~(\ref{Claw}), a three-step regularization scheme is devised to suppress the error and instabilities due to: (1)~\emph{compactness} of the double-layer potential~$\breve\bM\brbu$ in~(\ref{FOD1}), physically stemming from the presence of interface waves on~$\breve\Gamma$; (2)~tangential differentiation of the recovered FOD along \emph{inexact} fracture surface~$\breve\Gamma$, and (3)~presence (if any) of areas on $\breve\Gamma$ with \emph{near-zero} FOD values. These three steps of regularization are elucidated in the sequel. To aid the discussion, it is convenient to recall~(\ref{SVDM}) and to similarly introduce the SVD of matrix~$\breve{\text{\bf T}}$ discretizing the integro-differential operator~$\breve{\boldsymbol{T}}$ in~(\ref{TBIE}) as 
\beq\lb{SVDT}
\breve{\text{\bf T}} \;=\; \text{\bf U}_{\text{\tiny T}} \boldsymbol{\Lambda}_{\text{\tiny T}} \text{\bf V}_{\text{\tiny T}}^*,
\vspace{-1 mm}
\eeq
where $\text{\bf U}_{\text{\tiny T}}$ (resp. $\text{\bf V}_{\!\text{\tiny T}}$) collects the left (resp. right) eigenvectors of~$\breve{\bf T}$, and $\boldsymbol{\Lambda}_{\text{\tiny T}}$ contains its singular values.  

%--------------------------------------------------------------------------------------------------%
\subsubsection{Suppression of interface waves on~$\breve\Gamma$.}
%--------------------------------------------------------------------------------------------------%
Recalling~(\ref{FOD2}) as the basis for FOD recovery and assuming the availability of sensory data for multiple incident fields, ($p=1,\dots P$), the idea is to \emph{synthetically recombine} the available scattered-field data $\tilde{\text{\bf u}}^{\infty,p}$ in order to minimize the participation of interface waves in the associated FOD profile $\llbracket \breve{\text{\bf u}} \rrbracket$. To do so, recall that ${\bf M}$ is excitation-independent and let~$\text{\bf U}_{\text{\tiny M,q}}$ ($\text{q}\!=\!1,2,\ldots Q$) denote the \emph{left eigenvectors} affiliated with the~$Q$ \emph{smallest singular values} of~$\breve{\text{\bf M}}$ that are suppressed by the regularization of~(\ref{FOD2}). On conveniently rewriting the latter equation for the $p$th incident field as $\boldsymbol{\Lambda}_{\text{\tiny M}} \text{\bf V}_{\text{\tiny M}}^*\llbracket \breve{\text{\bf u}} \rrbracket = \text{\bf U}_{\text{\tiny M}}^* \tilde{\text{\bf u}}^{\infty,p}$, the idea is to solve 
\beq
\lb{FOD3}
\text{\bf U}_{\text{\tiny M,q}}^* \Big[\sum_{p=1}^{P} g_p \, \tilde{\text{\bf u}}^{\infty,p} \Big] ~=~0, \qquad q=\overline{1,Q}, \vspace*{-1mm}
\eeq
for the (synthetic) source density~$\bg=(g_{\text{\tiny 1}},g_{\text{\tiny 2}},\dots, g_{\text{\tiny \it P}})$, designed to \emph{minimize} the participation of (unrecoverable) interface waves in the FOD on~$\breve\Gamma$ -- without prior knowledge of~$\breve\bK(\bxi)$. In this way, the solution error due to regularization of~(\ref{FOD2}) via e.g. truncated SVD is suppressed, while maintaining the stability of the solution by deploying the sensory data $\tilde{\text{\bf u}}^{\infty}=\sum_{p=1}^{P}g_p \tilde{\text{\bf u}}^{\infty,p}$ as the basis for solving~(\ref{FOD2}). In this setting, one proceeds with solving~(\ref{Claw}) for the specific stiffness profile~$\breve{\text{\bf k}}$ by setting $\breve{\text{\bf t}}^i = \sum_{p=1}^P g_p \, \breve{\text{\bf t}}^{i,p}$ where $\breve{\text{\bf t}}\ff_p\,^{i,p}$ is the free-field traction on $\breve\Gamma$ due to $p$th incident field. 

\begin{rem}
Assuming $P\!<\!Q$ in solving (\ref{FOD3}), the ``\emph{strength}''~$g_{p}$ of each incident field can be computed via least squares. When $P\!>\!Q$, on the other hand, one can take $Q$ out of $P$ available sets of sensory data to construct an even-determined problem.     
\end{rem}

%--------------------------------------------------------------------------------------------------%
\subsubsection{Regularization of~$\,\breve{\text{\bf T}}$.}
%--------------------------------------------------------------------------------------------------%

In general the linear operator $\breve{\boldsymbol{T}}\!:\tilde{H}^{1/2}(\breve\Gamma)\to H^{-1/2}(\breve\Gamma)$ introduced in~(\ref{TBIE}), that is discretized as~$\breve{\text{\bf T}}$, is constructed on the basis of the recovered fracture support~$\breve\Gamma$ and entails tangential differentiation according to~(\ref{dod}). As a result, $\breve{\text{\bf T}}$ may contain large singular values due to the approximate nature of $\breve\Gamma$ and its roughness, leading to the amplification of small errors in ${\llbracket \breve{\text{\bf u}} \rrbracket}$ while computing the right-hand-side of (\ref{Claw}). To ensure a stable solution in terms of~$\breve{\text{\bf k}}$, consider the minimal subspace of~$\mathbb{R}^{3N^\text{nds}}$, spanned by the first~$N$ \emph{right eigenvectors} $\text{\bf V}_{\!\text{\tiny T,n}}$ ($\text{n}\!=\!1,2,\ldots N$) of~$\breve{\text{\bf T}}$, that contributes to the construction of ${\llbracket \breve{\text{\bf u}} \rrbracket}$ below a designated error $\delta$ (a measure of error in FOD reconstruction), namely
\beq\lb{regT}
\bigg\| {{\llbracket \breve{\text{\bf u}} \rrbracket} - \sum_{\text{n}=1}^{N} ({\llbracket \breve{\text{\bf u}} \rrbracket} \sip \text{\bf V}_{\text{\tiny T,n}}) \text{\bf V}_{\text{\tiny T,n}}}\bigg\| \,\,<\, \delta.
\eeq
In this setting, the first term on the right-hand side in~(\ref{Claw}) can be regularized as
\beq\lb{regT2}
\breve{\text{\bf T}} \, \llbracket \breve{\text{\bf u}} \rrbracket \,\simeq\, 
\sum_{\text{n}=1}^{N} \big( \breve{\bf T} \llbracket \breve{\bf u} \rrbracket \sip \text{\bf U}_{\text{\tiny T,n}} \big) \text{\bf U}_{\text{\tiny T,n}},
\eeq
where $\text{\bf U}_{\text{\tiny T,n}}$ ($\text{n}\!=\!1,2,\ldots N$) are the left eigenvectors of $\breve{\text{\bf T}}$. 

%--------------------------------------------------------------------------------------------------%
\subsubsection{Treatment of vanishing FOD on~$\breve\Gamma$.}\lb{VFOD}
%--------------------------------------------------------------------------------------------------%

Solving~(\ref{Claw}) for the interfacial stiffness profile $\breve{\text{\bf k}}$, one may find that the coefficient matrix $\breve{\text{\bf A}}_{\llbracket \breve{\text{\bf u}} \rrbracket}$ in~(\ref{Claw2}) or $\breve{\text{\bf B}}_{\llbracket \breve{\text{\bf u}} \rrbracket}$ in~(\ref{Claw3}) is singular. To interpret this situation, one may recall~(\ref{Claw2}) and observe that the near-zero eigenvalues of $\breve{\text{\bf A}}_{\llbracket \breve{\text{\bf u}} \rrbracket}$ are affiliated with those collocation points~$\bxi\in\breve\Gamma$ where ${\llbracket \breve{\text{\bf u}} \rrbracket}(\bxi) \rightarrow \bzero$. One such example can be constructed by considering a penny-shaped fracture with locally orthotropic $\bK(\bxi)$, subjected to a pair of \emph{antiplane shear} incident waves~(\ref{plwa}) whose directions of incidence are symmetric with respect to the fracture plane, and whose polarization is parallel to the fracture plane; in this case the FOD can be shown to vanish along the entire fracture. To ensure a stable solution in situations where the reconstructed FOD vanishes over a subset of germane collocation points, (\ref{Claw}) can be solved via suitable regularization e.g.~by invoking Tikhonov regularization or truncated SVD. In this case, however, additional incident fields may be needed to quantitatively resolve~$\breve\bK(\bxi)$ at those locations. 

%--------------------------------------------------------------------------------------------------%
\section{Numerical implementation and results}\lb{Numexpp}
%--------------------------------------------------------------------------------------------------%

In light of the existing elastodynamic simulations~\cite{Fatemeh2015, Fatemeh2016} -- specifically focused on the geometric reconstruction of heterogeneous fractures via TS and GLSM, this section examines the effectiveness of the proposed 3-step approach (see Fig.~\ref{inversion3}) with a particular emphasis on the inversion of the specific stiffness profile. 
\begin{figure}[!bp]
\begin{center}
\includegraphics[width=0.56\linewidth]{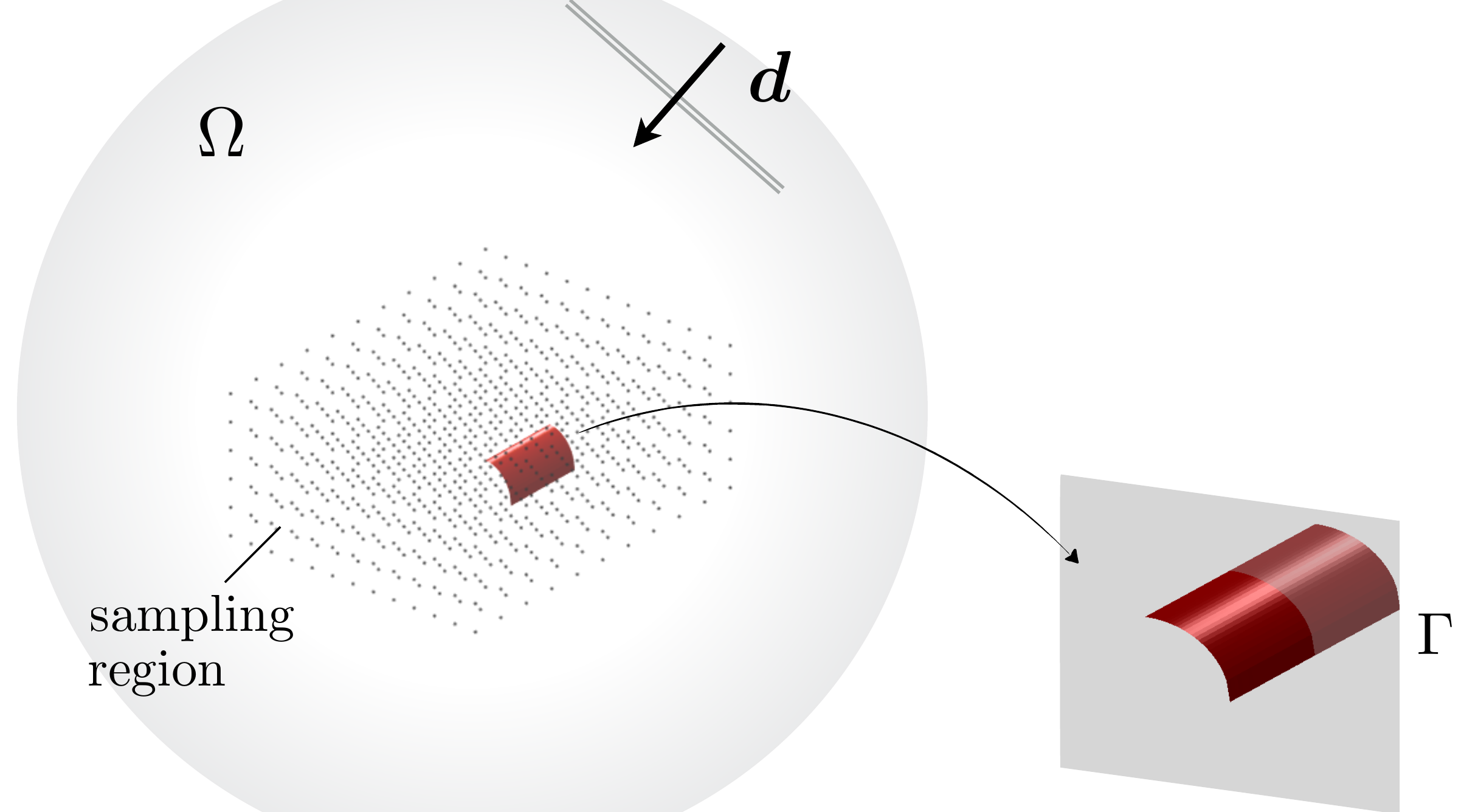}
\end{center} \vspace{-2mm}
\caption{\small{Elastodynamic sensing configuration featuring a heterogeneous cylindrical fracture.}}
\label{SC}
\end{figure}
The testing configuration is illustrated in Fig.~\ref{SC}, where a cylindrical fracture $\Gamma$ of width $L = 0.7$, arclength $\ell=0.55$, and radius $R = 0.35$ is embedded in a linear, isotropic and homogeneous elastic solid with shear and compressional wave speeds $c_s = 1$ and $c_p = 2.08$, respectively. As shown in Fig.~\ref{Int_true}, the fracture is endowed with two sets of (diagonal and orthotropic) interfacial stiffness profiles~$\bK(\bxi)$, namely: (1)~a piecewise-constant ``\emph{zebra}'' distribution in both shear and normal directions, and (2)~a ``\emph{cheetah}'' pattern defined in the fracture's local coordinates. By making reference to the local orthonormal basis $(\bn,\be_1,\be_2)$ on~$\Gamma$, one may conveniently write the specific stiffness matrix for both configurations as 
\beq\lb{zech}
\bK(\bxi) ~=~ \kappa_n(\bxi)\exs \bn\otimes\bn \;+\; \kappa_s(\bxi)\exs \big\{\be_1\otimes\be_1 \,+\, \be_2\otimes\be_2\}, \qquad \bxi\in\Gamma.
\eeq
As can be seen from Fig.~\ref{Int_true}, the pair $(\kappa_n,\kappa_s)$ is assumed to be \emph{complex-valued}, signifying energy dissipation (due to e.g.~friction) at the interface. The fracture is illuminated by a set of incident plane waves~\eqref{plwa}, taking~$k_s$ such that the ratio between the shear wavelength and the arclength of $\Gamma$ is $\lambda_s/\ell = 0.7$. Thus-induced scattered field is then measured in terms of the far-field pattern, $\tilde{\bf u}^\infty$, given by~\eqref{vinf}--\eqref{far-field}. The spatial density of sensory data, for both illumination and sensing purposes, is given by $N_\theta \!\times\! N_\phi = 25 \!\times\! 12$ directions given by the polar ($\theta_j,\, j\!=\!1,\ldots N_\theta$) and azimuthal~($\phi_k, \,k\!=\!1,\ldots N_\phi$) angle values spanning the unit sphere~$\Omega$. In what follows, all forward simulations are performed by an elastodynamic boundary integral equation code~\cite{Pak1999,Fatemeh2015} and polluted by 5\% random noise in terms of the observed far-field pattern.

\emph{Geometric reconstruction.} For each interface scenario the fracture support $\breve\Gamma$ is recovered, by way of the GLSM framework described in Section~\ref{GLSM}, from the knowledge of the far-field operator~\eqref{ffo0} as shown in Fig.~\ref{georec} (see~\cite{Fatemeh2016} for details).

\begin{figure}[!tp]
\begin{center}
\includegraphics[width=0.72\linewidth]{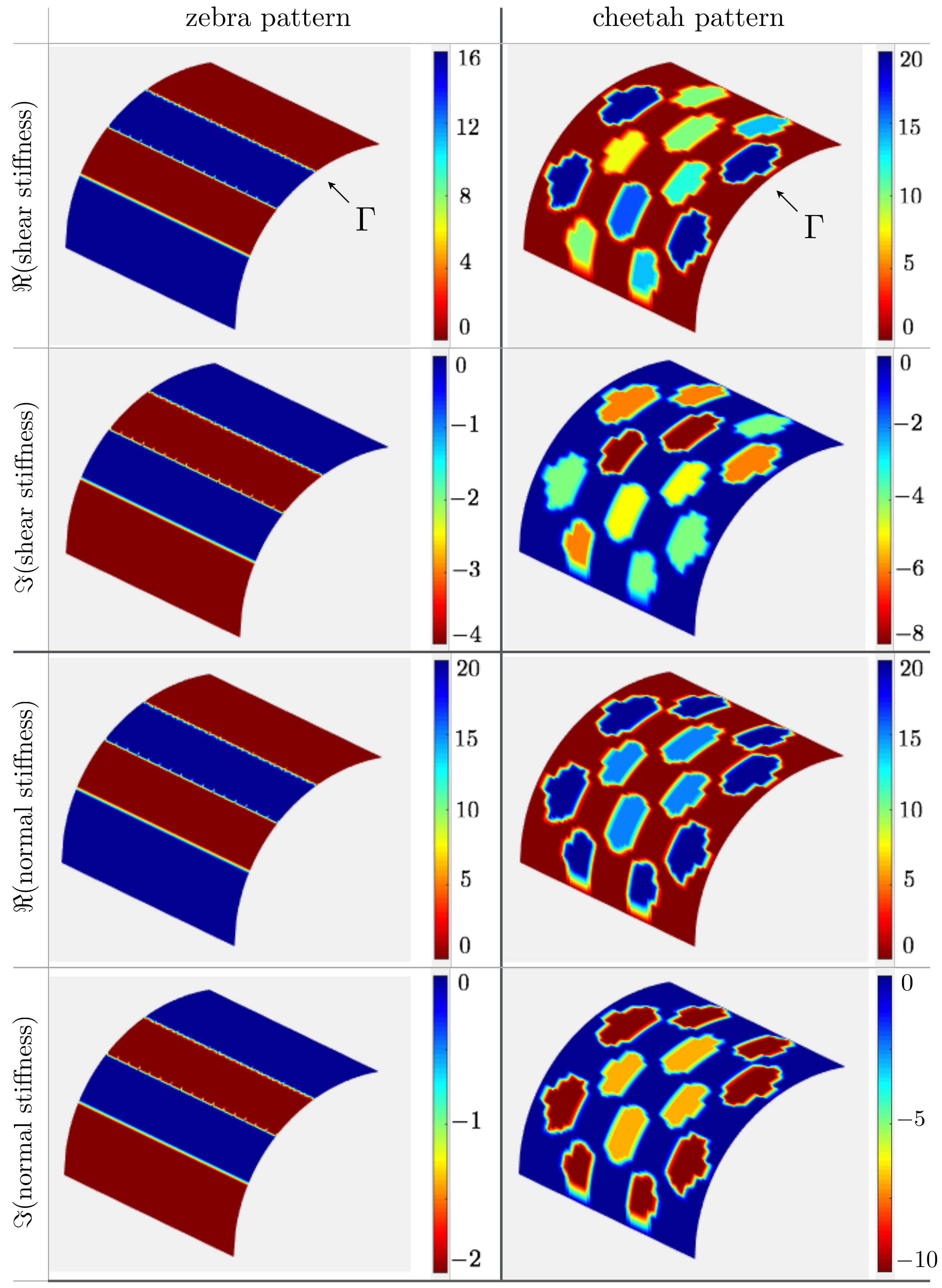}
\end{center} \vspace{-4 mm}
\caption{\small{Two interface scenarios for the example cylindrical fracture: \emph{zebra} (left) and \emph{cheetah} (right) patterns representing the distribution of complex-valued specific stiffness coefficients~$k_s$ and~$k_n$ in~\eqref{zech}.}}
\label{Int_true}
\end{figure}
\begin{figure}[!h]
\begin{center}
\includegraphics[width=0.65\linewidth]{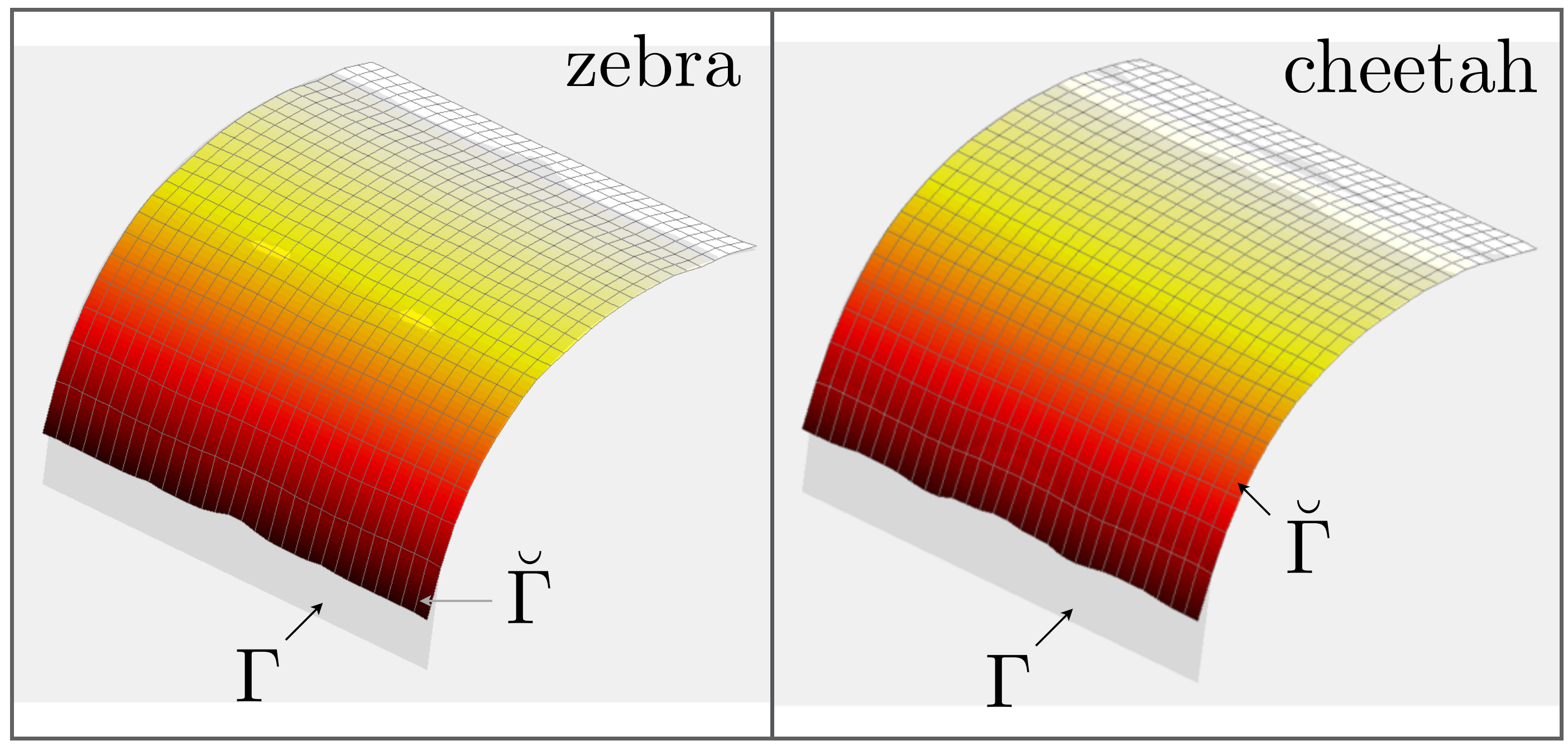}
\end{center} \vspace*{-4mm}
\caption{\small{Geometric reconstruction by way of~GLSM: ``true" geometry $\Gamma$ vs.~the reconstructed fracture support $\breve\Gamma$, obtained for the zebra and cheetah interface scenarios.}}
\label{georec}\vspace*{-2mm}
\end{figure}
 
\begin{figure}[!bp]
\begin{center}
\includegraphics[width=0.75\linewidth]{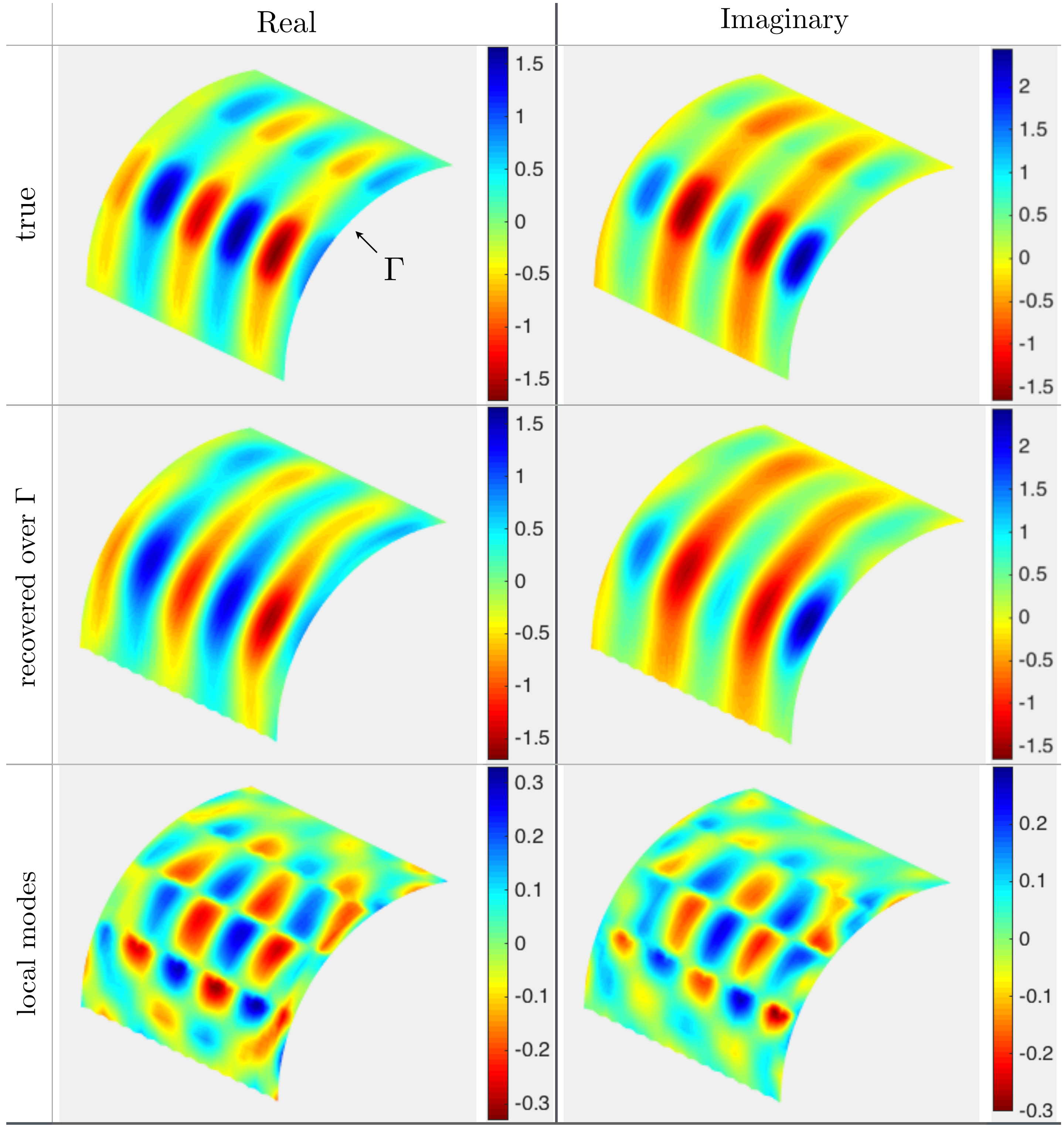}
\end{center} \vspace{-4mm}
\caption{\small{FOD recovery assuming prior knowledge of the ``true'' fracture geometry~$\Gamma$: shear component of $\llbracket\breve{\text{\bf u}}\rrbracket$ (along the width of $\Gamma$), \emph{zebra} interface scenario, single incident wave. The ``true'' FOD, $\llbracket{\bf u}\rrbracket$ (top row) consists of (i) the trace~$\llbracket\breve{\bf u}\rrbracket$ of the propagating waves (middle row) -- whose non-trivial fingerprint in the far-field data allows for their robust reconstruction, and (ii) the residual part $\llbracket{\bf u}\rrbracket-\llbracket\breve{\bf u}\rrbracket$ (bottom row), given by the trace of local modes i.e.~interface waves -- whose vanishing signature in the far-field data prevents their recovery.}}
\label{FOD_int_wave}
\end{figure}

\begin{figure}[!tp]
\begin{center}
\includegraphics[width=0.92\linewidth]{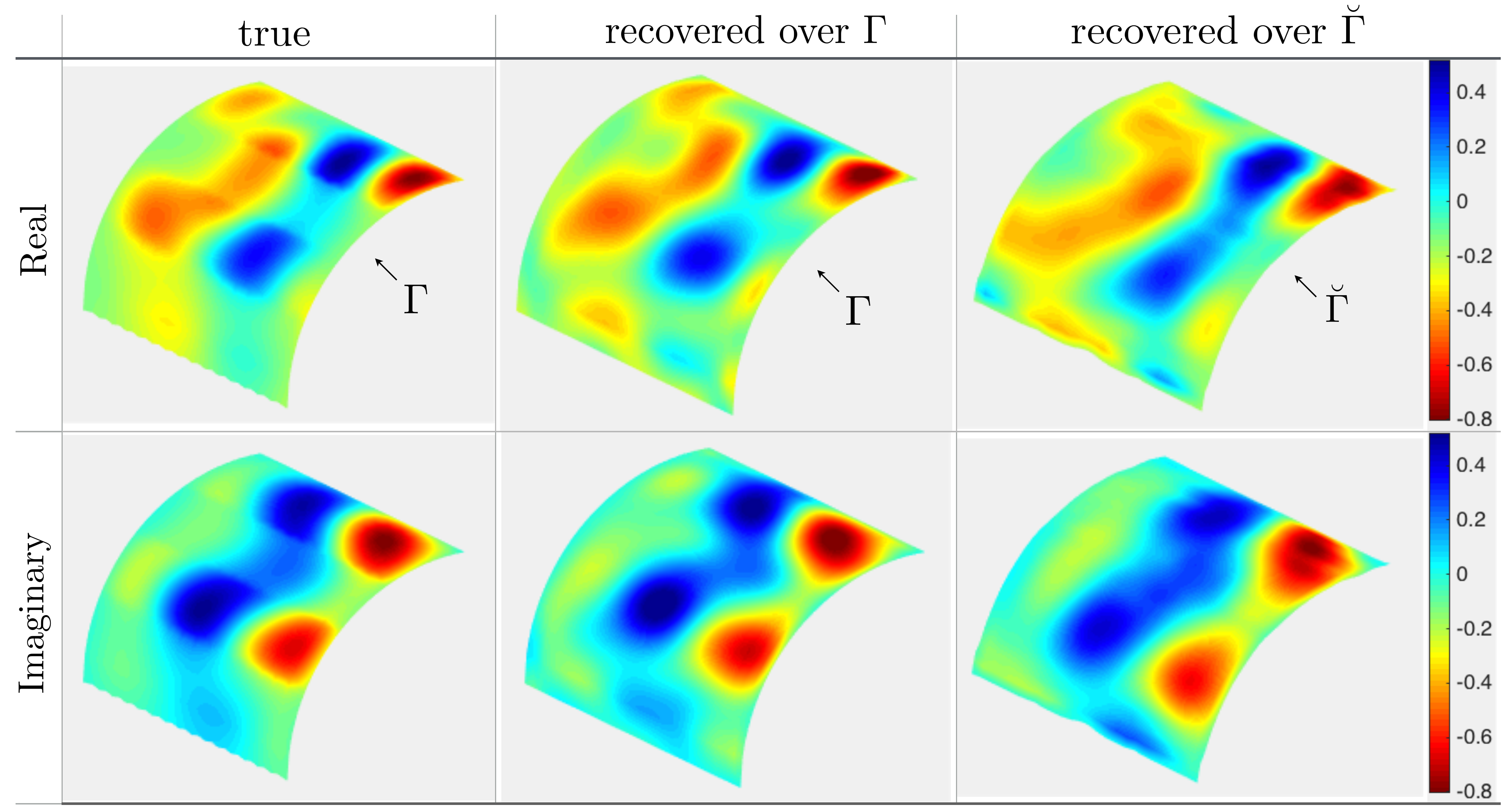}
\end{center} \vspace{-5 mm}
\caption{\small{FOD recovery from the far-field data collected for a \emph{single} incident P-wave, normal component of  $\llbracket\breve{\text{\bf u}}\rrbracket$, \emph{zebra} interface scenario: ``true'' FOD on $\Gamma$ as given by the forward model (left); recovered FOD over the true fracture geometry $\Gamma$ (middle); recovered FOD over the reconstructed fracture geometry $\breve\Gamma$ (right).}}
\label{FOD_rec_v2}
\vspace{9mm}
\begin{center}
\includegraphics[width=0.78\linewidth]{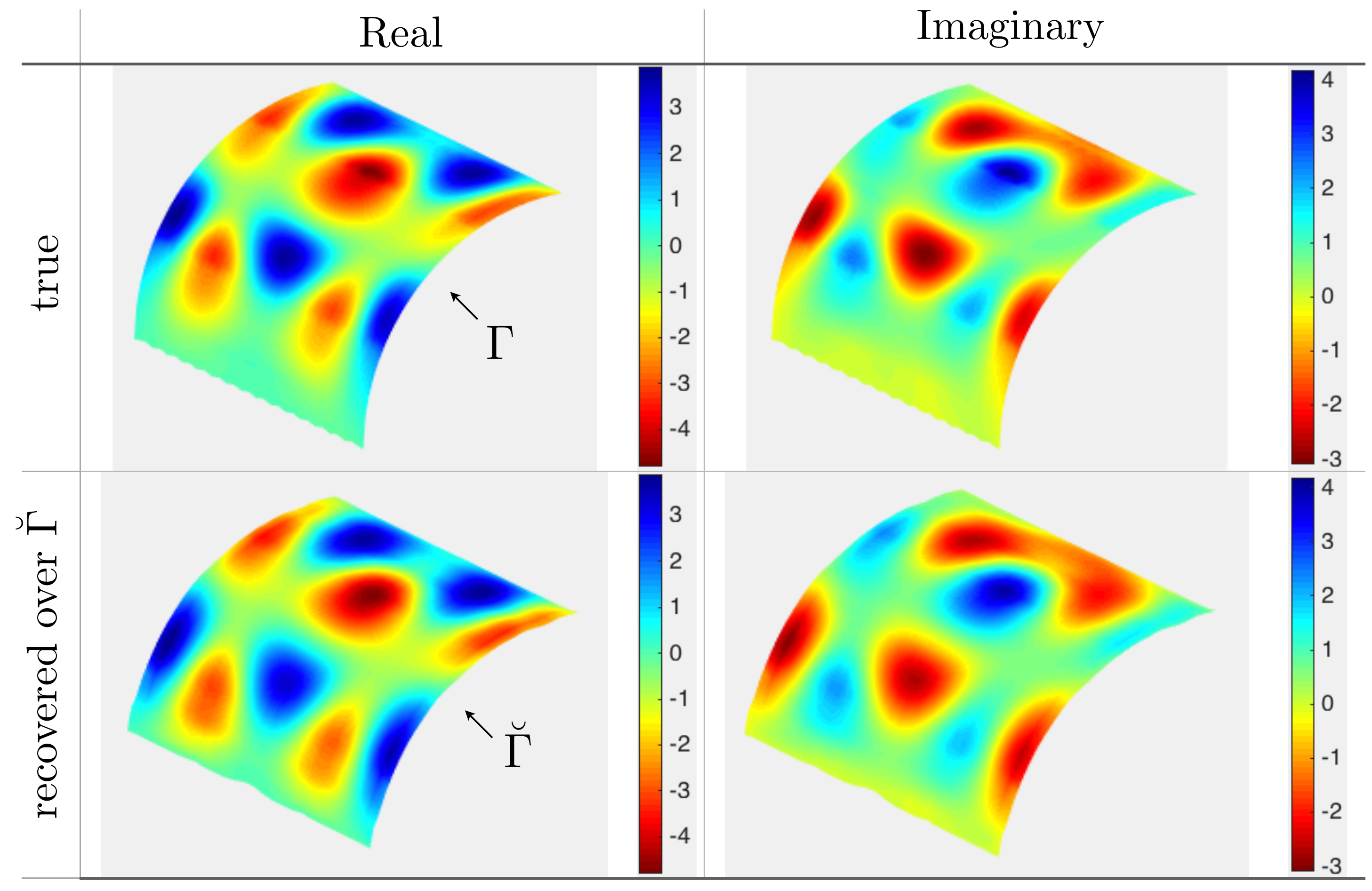}
\end{center} \vspace{-5mm}
\caption{\small{FOD recovery from the far-field data collected for \emph{multiple} incident waves, normal component of  $\llbracket\breve{\text{\bf u}}\rrbracket$, \emph{zebra} interface scenario: ``true'' FOD on $\Gamma$ (top) induced by a set of incident plane waves whose density is obtained via~(\ref{FOD3}) vs. recovered FOD over the reconstructed fracture surface $\breve\Gamma$ (bottom).}}
\label{FOD_rec_fix_v2}
\end{figure}

\emph{FOD recovery.}~Given $\breve\Gamma$, one may construct a ``double-layer potential'' operator~$\breve{\text{\bf M}}$ in~\eqref{FOD2} -- mapping the sought FOD, ${\llbracket \breve{\bf u} \rrbracket}$, to the sensory data $\tilde{\bf u}^\infty$ by way of~\eqref{FOD1}. As shown by Lemma~\ref{compact}, \eqref{FOD2} represents an ill-posed problem due to the compactness of~$\breve{\bM}$ which is reflected by the appearance of interfacial waves~\cite{Pyr1987} on~$\breve{\Gamma}$ -- which are known to decay exponentially away from the fracture surface. This is illustrated in Fig.~\ref{FOD_int_wave} where the ``true'' FOD over $\Gamma$ -- obtained by simulating the forward problem~(\ref{GE}) due to a single incident S-wave, is decomposed into two components, namely: (i)~the \emph{retrievable} part~${\llbracket \breve{\text{\bf u}} \rrbracket}$ computed by solving~(\ref{FOD2}) over $\Gamma$ via Tikhonov regularization endowed with the Morozov discrepancy principle ($5\%$ random noise), and (ii)~the \emph{residual} part~${\llbracket{\bf u}\rrbracket}-{\llbracket \breve{\text{\bf u}} \rrbracket}$, comprised of interfacial waves, obtained by subtracting the reconstructed displacement jump from the ``true'' FOD. Note that $\breve{\text{\bf M}}$ used to recover FOD in Fig.~\ref{FOD_int_wave} is intentionally constructed on the basis of the \emph{``true''} fracture geometry $\Gamma$, so that the computed residual part is not polluted by errors due to geometric fracture reconstruction. 
With such remark, Fig.~\ref{FOD_rec_v2} compares the ``true'' FOD over $\Gamma$, ${\llbracket{\bf u}\rrbracket}$ -- induced by the interaction of a single incident P-wave with the ``zebra" interface on $\Gamma$, with the recovered displacement jumps, ${\llbracket\breve{\bf u}\rrbracket}$ according to~(\ref{FOD2}), over~$\Gamma$ and~$\breve\Gamma$. For a robust inversion of the interfacial stiffness, however, one should reconstruct the FOD profile from~$\tbu^\infty$ by making use of the first regularization step in Section~\ref{sstif} -- where the observed wavefield data from multiple incident fields are \emph{synthetically recombined} to minimize, and possibly eliminate, the participation of interface waves in the FOD on~$\breve{\Gamma}$. With reference to~(\ref{FOD3}), this is accomplished by selecting~$Q$ as the number of singular values of $\breve{\text{\bf M}}$ that are smaller than $15\%$ of its largest singular value. The result is is shown in Fig.~\ref{FOD_rec_fix_v2}, which compares the ``true'' FOD over $\Gamma$ -- due to recombined incident fields -- with the corresponding reconstruction  ${\llbracket \breve{\text{\bf u}} \rrbracket}$ over $\breve\Gamma$. 

\begin{figure}[!tp]
\begin{center}
\includegraphics[width=0.74\linewidth]{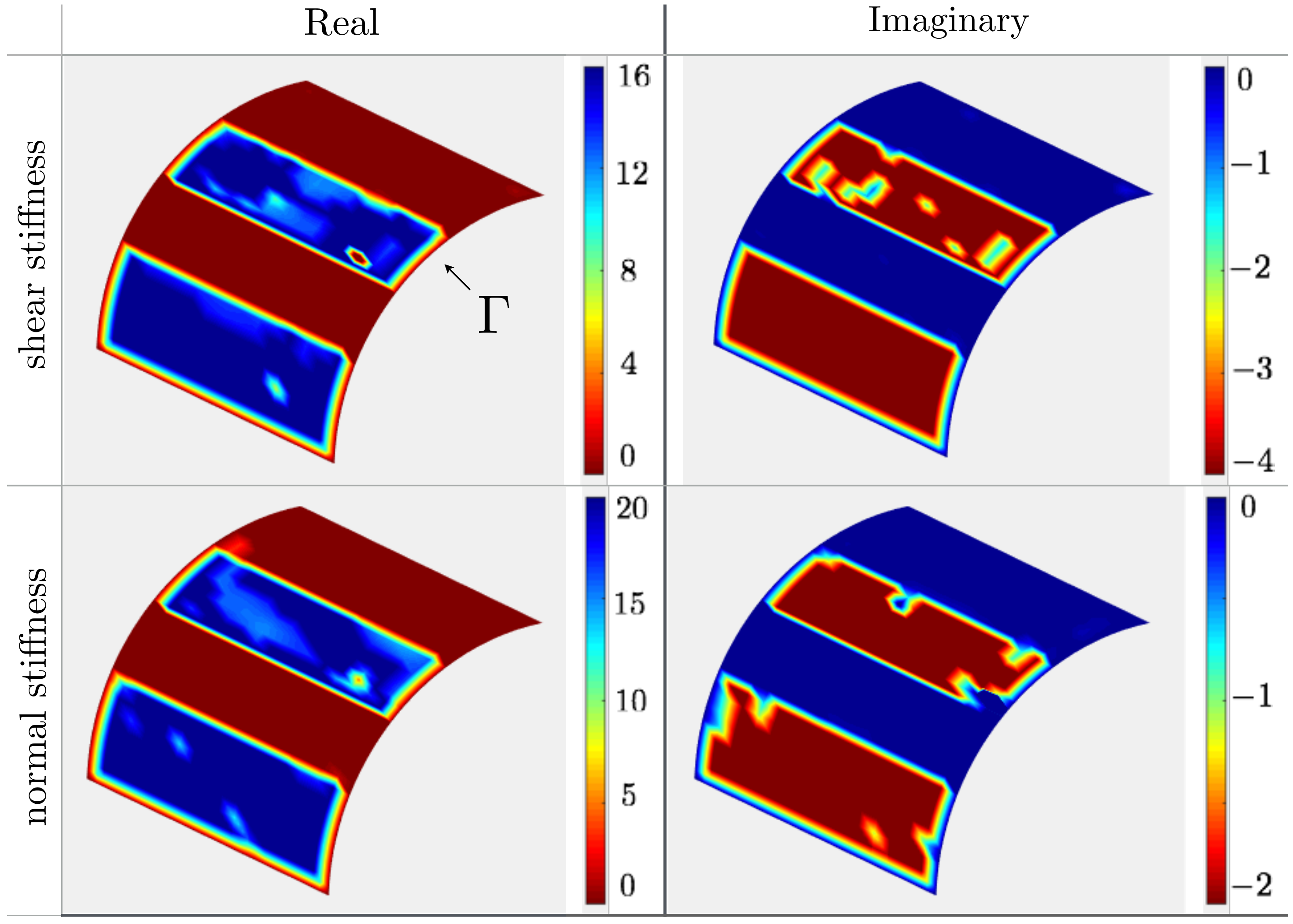}
\end{center} \vspace{-5mm}
\caption{\small{Recovered distribution of the specific stiffness, assuming prior knowledge of the ``true'' fracture geometry $\Gamma$, from the far-field patterns collected for multiple incident waves: normal and shear components of $\breve{\bK}(\bxi)$, \emph{zebra} interface scenario.}}
\label{Int_rec_true}
\end{figure} 

\begin{figure}[!tp]
\begin{center}
\includegraphics[width=0.70\linewidth]{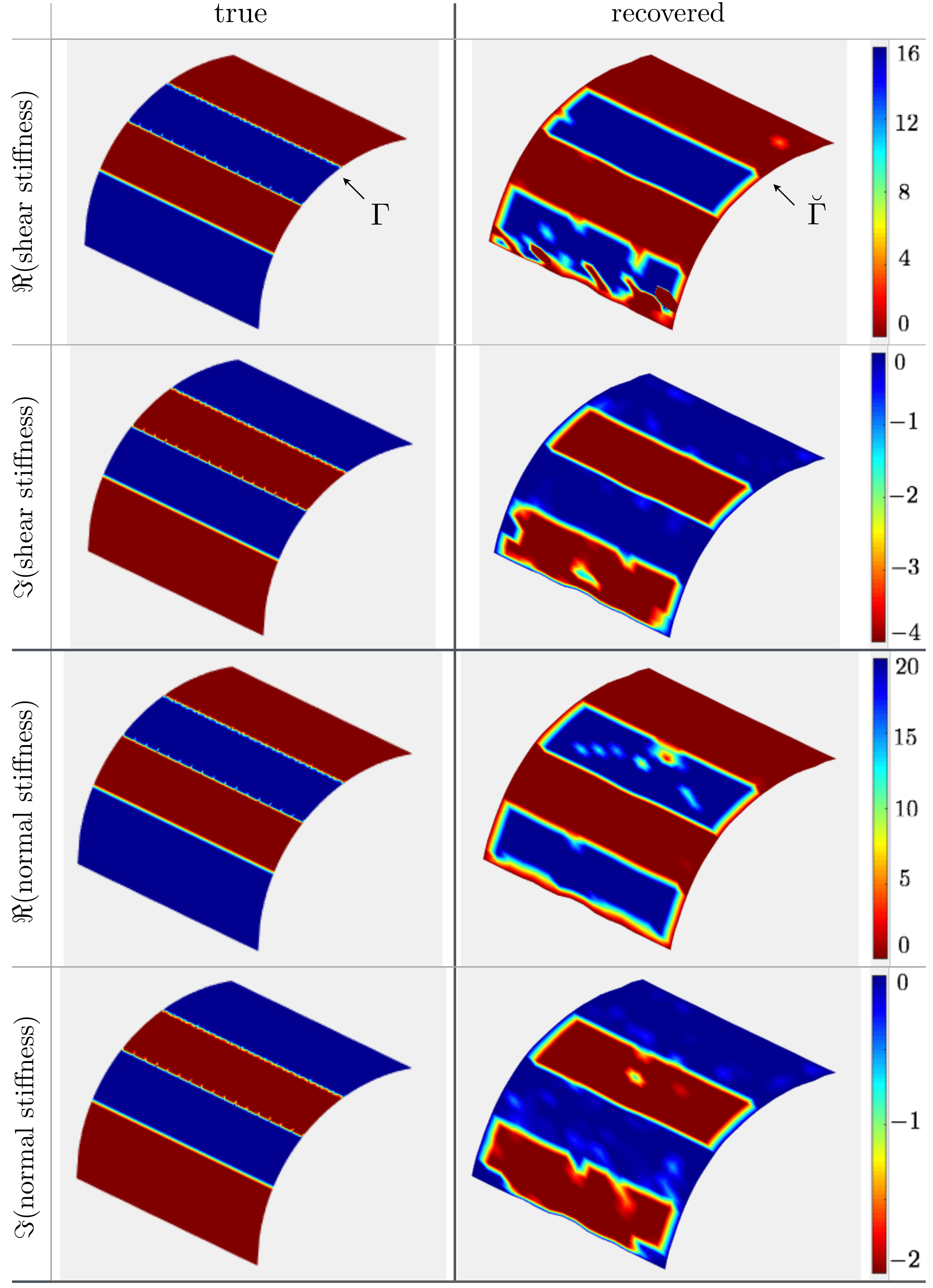}
\end{center} \vspace{-5mm}
\caption{\small{Recovered distribution of the specific stiffness, over the reconstructed fracture surface~$\breve\Gamma$, from the far-field patterns collected for multiple incident waves: ``true'' \emph{zebra} pattern $\bK(\bxi)$ over~$\Gamma$ (left) vs. the recovered distribution $\breve{\bK}(\bxi)$ over $\breve\Gamma$ (right).}}
\label{Int_rec_rec_Zebra}
\end{figure}

\begin{figure}[!h]
\begin{center}
\includegraphics[width=0.72\linewidth]{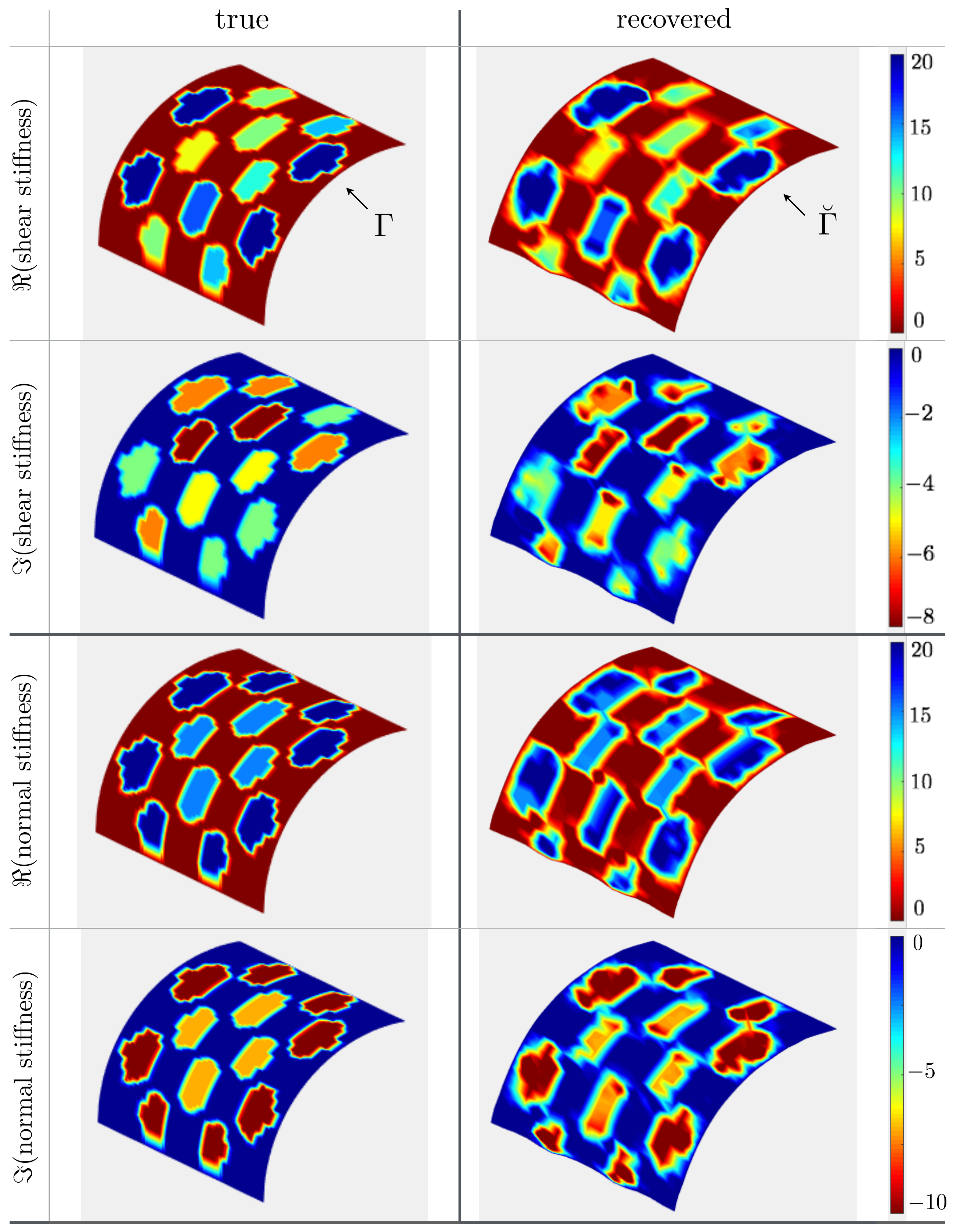}
\end{center} \vspace{-5mm}
\caption{\small{Recovered distribution of the specific stiffness, over the reconstructed fracture surface~$\breve\Gamma$, from the far-field patterns collected for multiple incident waves: ``true'' \emph{cheetah} pattern $\bK(\bxi)$ over~$\Gamma$ (left) vs. the recovered distribution $\breve{\bK}(\bxi)$ over $\breve\Gamma$ (right).}}
\label{Int_rec_rec_cheetah}
\end{figure}

\emph{Reconstruction of the specific stiffness.} By virtue of~\eqref{zech} and the developments in Section~\ref{sstif}, one finds the regularized statement of~\eqref{Claw} to read 
\beq\lb{Kf}
\breve{\text{\bf A}}_{\llbracket \breve{\text{\bf u}} \rrbracket} \exs \breve{\bf k} ~=~  
\sum_{\text{n}=1}^{N} \big(\breve{\text{\bf T}} \llbracket \breve{\text{\bf u}} \rrbracket  \sip \text{\bf U}_{\text{\tiny T,n}}\big) \text{\bf U}_{\text{\tiny T,n}}  \,+\, \breve{\text{\bf t}}^{i}.
\eeq           
Here, $\breve{\text{\bf A}}_{\llbracket \breve{\text{\bf u}} \rrbracket}$ a diagonal coefficient matrix given by the reconstructed FOD at collocation points on~$\breve{\Gamma}$; $\breve{\bf t}^i$ collects the incident-field tractions on $\breve\Gamma$, and $\text{\bf U}_{\text{\tiny T,n}}$ ($\text{n}\!=\!1,2,\ldots N$) are the left eigenvectors of $\breve{\text{\bf T}}$, truncated as in~(\ref{regT}) with~$\delta = 0.001$ to minimize the effect of imperfect geometric reconstruction. As examined in Section~\ref{VFOD}, (\ref{Kf}) is solved via Tikhonov regularization and Morozov discrepancy principle applied to the noise level of $5\%$. In this setting, Fig.~\ref{Int_rec_true} shows the reconstructed \emph{zebra} pattern over $\Gamma$ -- assuming prior knowledge of the ``true" fracture geometry $\Gamma$, while Fig.~\ref{Int_rec_rec_Zebra} (resp.~Fig.~\ref{Int_rec_rec_cheetah}) compares the ``true" \emph{zebra} (resp.~\emph{cheetah}) $\bK$-distributions over $\Gamma$ with their recovered $\breve\bK$-counterparts on $\breve\Gamma$. As can be seen from both displays, the fidelity of specific stiffness reconstruction is rather remarkable given (i) no prior information about the fracture geometry nor its contact condition, and (ii) multiple steps of regularization used to accomplish the sequential geometric reconstruction and interfacial characterization. 

%--------------------------------------------------------------------------------------------------%
\section{Acknowledgments}
%--------------------------------------------------------------------------------------------------% 

The second author kindly acknowledges the support provided by the National Science Foundation (Grant \#1536110) and the University of Minnesota Supercomputing Institute.

%--------------------------------------------------------------------------------------------------%
\section{Summary}
%--------------------------------------------------------------------------------------------------%

In this study, a three-step inverse solution is proposed for the geometric reconstruction and interfacial characterization of heterogeneous fractures by seismic waves. In the approach the fracture surface is (a)~allowed to be unconnected, and (b)~endowed with spatially-dependent (linear) contact condition via a $3\times3$ complex-valued matrix of specific stiffness coefficients, which is capable of representing the effects of dilatancy, energy dissipation, and anisotropy in the fracture's contact law. In the first stage of the sensing scheme, the fracture surface is reconstructed with the aid of recently established techniques for elastic waveform tomography of heterogeneous fractures (namely the methods of topological sensitivity and generalized linear sampling) which operate irrespective of the fracture's (unknown) contact condition. With such result in place, the fracture's opening displacement (FOD) is, for given sensory data, recovered over the reconstructed fracture surface~$\breve\Gamma$ via a double-layer potential representation mapping the FOD to the scattered field in the exterior domain. In the third and last step of the proposed scheme the specific stiffness profile (as given by its normal, shear, and mixed-mode components) -- is recovered from the knowledge of~FOD and~$\breve\Gamma$ by making use of the traction boundary integral equation written for~$\breve\Gamma$. To help construct a robust inverse solution, a three-step regularization scheme is also devised to minimize the error due to: (i)~compactness of the double-layer potential map used to recover the FOD, physically manifested by the presence of interfacial waves on~$\breve\Gamma$; (ii)~inexact geometric reconstruction, and (iii)~possible presence of areas on $\breve\Gamma$ with near-zero FOD values. The numerical results, assuming no prior knowledge of the fracture geometry nor its contact condition, suggest a remarkable performance of the staggered inverse solution and its potential toward quantitative recovery of both the elastic and dissipative contact characteristics of heterogeneous fractures.  
 
\newpage

%=======================================================================================================%
\bibliography{inverse,crackbib}
%=======================================================================================================%

%*******************************************************************************************************%
\end{document}